\documentclass[11pt]{article}

\usepackage{amsmath}
\usepackage{amssymb}
\usepackage{amsthm}
\usepackage{latexsym}
\usepackage{color}
\usepackage{graphicx}
\usepackage{color}

\DeclareSymbolFont{calletters}{OMS}{cmsy}{m}{n}
\DeclareSymbolFontAlphabet{\mathcal}{calletters}

\newtheorem{Theorem}{Theorem}[part]
\newtheorem{Definition}{Definition}[part]

\newtheorem{Assumption}{Assumption}[part]

\newtheorem{Remark}{Remark}[part]

\newtheorem{Properties}{Properties}

\makeatletter \@addtoreset{equation}{section}

\@addtoreset{Definition}{section}

\@addtoreset{Theorem}{section}

\@addtoreset{Proposition}{section}

\@addtoreset{Properties}{section}

\@addtoreset{Property}{section}

\@addtoreset{Assumption}{section}

\@addtoreset{Assertion}{section}

\@addtoreset{Corollary}{section}

\@addtoreset{Lemma}{section}

\@addtoreset{Remark}{section}

\@addtoreset{Example}{section}

\def \proof{{\noindent \bf Proof. }}

\newcommand{\Pc}{\mathcal{P}}
\newcommand{\Pset}{\mathcal{P}_{H}}

\newcommand{\Pbb}{\mathbb{P}}
\newcommand{\Ebb}{\mathbb{E}}

\def \E{\mathbb{E}}
\def \F{\mathbb{F}}

\def \P{\mathbb{P}}
\def \Q{\mathbb{Q}}
\def \R{\mathbb{R}}
\def\Ac{{\cal A}}

\def\Ec{{\cal E}}

\def\Lc{{\cal L}}

\def \ep{\hbox{ }\hfill{ ${\cal t}$~\hspace{-5.1mm}~${\cal u}$   } }
\newcommand{\No}[1]{\left\|#1\right\|}     % Norm
\newcommand{\abs}[1]{\left|#1\right|}     % absolute value
\def\einf{{\rm ess \, inf}}
\def\esup{{\rm ess \, sup}}
\def\reff#1{{\rm(\ref{#1})}}

\def\log{{\rm log}}
\def\exp{{\rm exp}}

\def\sup{{\rm sup}}
\def\inf{{\rm inf}}

 \title{Robust Utility Maximization in Non-dominated Models with 2BSDEs\footnote{Research supported by the Chair {\it Financial Risks} of the {\it Risk Foundation} sponsored by Soci\'et\'e G\'en\'erale, the Chair {\it Derivatives of the Future} sponsored by the {F\'ed\'eration Bancaire Fran\c{c}aise}, and the Chair {\it Finance and Sustainable Development} sponsored by EDF and Calyon.}
}
\author{ Anis {\sc Matoussi}\footnote{CMAP, Ecole Polytechnique, Paris, and Universit\'e du Maine, Le Mans, anis.matoussi@univ-lemans.fr.}
			\and
			Dylan {\sc Possama\"i}\footnote{CMAP, Ecole Polytechnique, Paris, dylan.possamai@polytechnique.edu.}
      \and
      Chao {\sc Zhou}\footnote{CMAP, Ecole Polytechnique, Paris, chao.zhou@polytechnique.edu.}
      }
 \date{{\small{first submission  January 4, 2012\\
 revision version May 16, 2012}}}

\begin{document}

 \maketitle

 \begin{abstract}

The problem of robust utility maximization in an incomplete market with volatility uncertainty is considered, in the sense that the volatility of the market is only assumed to lie between two given bounds. The set of all possible models (probability measures) considered here is non-dominated. We propose studying this problem in the framework of second order backward stochastic differential equations (2BSDEs for short) with quadratic growth generators. We show for exponential, power and logarithmic utilities that the value function of the problem can be written as the initial value of a particular 2BSDE and prove existence of an optimal strategy. Finally several examples which shed more light on the problem and its links with the classical utility maximization one are provided. In particular, we show that in some cases, the upper bound of the volatility interval plays a central role, exactly as in the option pricing problem with uncertain volatility models of \cite{alp}.

\vspace{10mm}

\vspace{5mm}


\end{abstract}

\newpage

\section{Introduction}

One of the most prominent problems of mathematical finance literature is the so-called problem of utility maximization. It is a problem of optimal investment faced by an economic agent who has the opportunity to invest in a financial market consisting of a riskless asset and (for simplicity) one risky asset. Given a fixed investment horizon $T$, the aim of the agent is to find an optimal allocation between the two assets, so as to maximize his "welfare" at time $T$. Following the seminal work of Von Neumann and Morgenstern \cite{vnm}, where they assumed that the preference of the agent could be represented by a utility function $U$ and a given probability measure $\mathbb P$ reflecting his views, the now classical formulation of the problem consists in solving the optimization
$$V(x):=\underset{\pi\in\mathcal A}{\sup}\text{ }\mathbb E^\mathbb P\left[U(X_T^{x,\pi}-\xi)\right],$$
where $\mathcal A$ is the set of admissible strategies $\pi$ for the agent, $X_T^{x,\pi}$ is his wealth at time $T$ with initial capital $x$ and a trading strategy $\pi$, and $\xi$ is a terminal liability.

\vspace{0.8em}
In the sixties, Merton \cite{mer69} was the first to study and solve this problem in the particular case where the risky asset follows a Black-Scholes model, where there are no restrictions on the admissible strategies (that is to say in a complete market), where the utility function is of power type and where the liability is equal to $0$. The proof relies on classical techniques of stochastic control theory, since he manages to solve the Hamilton-Jacobi-Bellman PDE associated with the problem explicitly, and he then uses a verification argument. The problem in complete markets but with general utility functions was only solved in the eighties by Pliska \cite{plis}, using techniques from convex duality. Following these papers, a large trend of literature tried to weaken their assumptions, and notably the completeness hypothesis on the market, which was too restrictive and unrealistic from the point of view of applications. One possible direction of generalization is to impose constraints on the strategies of the investor. Following the first works of Cvitani\'c and Karatzas \cite{cvi} and Zariphopoulou \cite{zar}, where once more convex duality techniques were used, the beginning of the $21$st century saw the emergence of a link between this optimal investment problem and the theory of backward stochastic differential equations (BSDEs for short). These objects were first introduced by Bismut \cite{bis} in the linear case, then generalized by Pardoux and Peng \cite{pardpeng} to Lipschitz generators. On a filtered probability space $(\Omega,\mathcal F,\left\{\mathcal F_t\right\}_{0\leq t\leq T},\mathbb P)$ generated by a $\mathbb R^d$-valued Brownian motion $B$, a solution to a BSDE consists of a pair of progressively measurable processes $(Y,Z)$ such that
$$Y_t=\xi +\int_t^Tf_s(Y_s,Z_s)ds-\int_t^TZ_sdB_s,\text{ } t\in [0,T], \text{ }\mathbb P-a.s.$$
where $f$ (also called the driver) is a progressively measurable function and $\xi$ is a $\mathcal F_T$-measurable random variable.

\vspace{0.8em}
Then, El Karoui and Rouge \cite{ekr} considered the problem of indifference pricing with an exponential utility function (which is linked to the optimal investment problem) in the case where the strategies are constrained to stay in a given closed and convex set. They proved that the value function of the problem was related to the initial value of a BSDE with a driver of quadratic growth in the $Z$ part. Building upon these results, Hu, Imkeller and M\"uller \cite{him} generalized the approach to the case of logarithmic and power utilities with strategies constrained in a closed set.

\vspace{0.8em}
Another direction of generalization of the original Merton problem is related to the question of model uncertainty. Indeed, in all the above formulations, a probability measure $\mathbb P$ is fixed. It means that the investor knows the "historical" probability $\P$ that describes the dynamics of the state process. In reality, the investor may have some uncertainty on this probability, which means that there can be several objective probability measures to consider. The problem then becomes  a robust utility maximization and can be written as follows
$$V^{\xi}(x)\;:=\; \underset{\pi\in\mathcal{A}}{{\rm sup}}\ \underset{{\Q\in\Pc}}{{\rm inf}}\Ebb^{\Q}[U(X^{x,\pi}_T-\xi)],$$
where $\Pc$ is the set of all considered possible probability measures.

\vspace{0.8em}
In this case, the properties of the set $\mathcal P$ become crucial in order to solve the problem. The first results in the literature were limited to dominated sets. A set $\mathcal P$ is said to be dominated if every probability measure $\mathbb P\in\mathcal P $ is absolutely continuous with respect to some reference probability measure in $\mathcal P$. For instance, this is the case when considering drift uncertainty. In this framework, the problem was introduced by Gilboa and Schmeidler \cite{gs}. Anderson, Hansen and Sargent \cite{ahs} and Hansen et al. \cite{hstw} then introduced and discussed the basic problem of robust utility maximization, penalized by a relative entropy term of the model uncertainty $\Q \in \mathcal P$ with respect to a given reference probability measure $\P_0$. Inspired by these latter works, Bordigoni, Matoussi and Schweizer \cite{bms} solved the robust problem (the minimization part) in a more general semimartingale framework by using stochastic control techniques and proved that the solution was related to a quadratic semimartingale BSDE. Among others, results in the robust maximization problem were also  obtained by Gundel \cite{gun} , Schied and Wu \cite{sw} or Skiadas \cite{ski} in the case of continuous filtrations. The overall approach relies essentially on convex duality ideas.

\vspace{0.8em}
The situation becomes more intricate when the set $\mathcal P$ is no longer dominated, which happens when introducing volatility uncertainty, in the sense that the volatility process is only assumed to lie between two given bounds. Although the problem of option pricing under volatility uncertainty has been solved for a long time (see \cite{alp} and \cite{lyo} for instance), the problem of utility maximization was not addressed until recently by Denis and Kervarec \cite{dk} (see also \cite{tz} where it is analyzed under the framework of stochastic games and Hamilton-Jacobi-Bellman-Isaacs equations). In this article, they first establish a duality theory for robust utility maximization and then show that there is a least favorable probability measure and an optimal strategy. However, their utility function $U$ is supposed to be bounded %( thus not including exponential, power and logarithmic utility)
and to satisfy Inada conditions. Recently, in \cite{ej}, Epstein and Ji formulate a model of utility in a continuous-time framework that captures the decision-maker's concern with ambiguity or model uncertainty, even though they do not study the maximization problem of robust utility \emph{per se}. More recently, Tevzadze et al. \cite{tev2} studied a related robust utility maximization problem for exponential and power utility functions (and also for mean-square error criteria). We will compare our results and theirs in Section \ref{6} (see Remark \ref{tev}).

\vspace{0.8em}
The intuition at the core of our work is that, exactly as the problem of utility maximization under constraints was linked to BSDEs with quadratic growth, the problem of robust utility maximization under volatility uncertainty should be linked to some kind of backward equations. In fact, the right objects to consider in that case are the so-called second-order BSDEs (2BSDEs for short) which were introduced for the first time by Cheredito, Soner, Touzi and Victoir \cite{cstv}. However, they were not able to provide a complete theory of existence and uniqueness. Hence, a reformulation was proposed by Soner, Touzi and Zhang \cite{stz}, who provided a wellposedness theory for 2BSDEs under uniform Lipschitz conditions similar to those of Pardoux and Peng. Their key idea was to reinforce the condition that the 2BSDE must hold $\mathbb P-a.s.$ for every probability measure $\mathbb P$ in a non-dominated class of mutually singular measures (see Section \ref{section.1rob} for precise definitions). The theory being very recent, the literature remains rather limited. However, we refer the interested reader to Possama\"i \cite{pos} and Possama\"i and Zhou \cite{posmz1} who respectively extended these wellposedness results to generators with linear and quadratic growth. In our main result, we show that in incomplete markets with volatility uncertainty, the solution of the robust utility maximization problem (for exponential, logarithmic and power utilities) is related to the initial value of a particular 2BSDE with quadratic growth, as suggested by the intuition. We also emphasize a specificity in our approach when it comes to the sets of admissible strategies considered. Usually, when dealing with this type of problems (see for instance \cite{ekr} and \cite{him}), an exponential uniform integrability assumption is made on the trading strategies. Our approach relies instead on integrability assumptions of BMO type on the trading strategies. The mathematical justifications are detailed in Remarks \ref{bilou} and \ref{bmoK}. However, there is also a financial interpretation. Indeed, as explained in \cite{frei} which adopts the same type of BMO framework, this assumption corresponds to a situation where the market price of risk is assumed to be BMO. Exactly as in the case of a bounded market price of risk, it implies that the minimum martingale measure is a true probability measure, and therefore that the market is without arbitrage, in the sense of No Free Lunch with Vanishing Risk.

\vspace{0.8em}
The rest of the paper is organized as follows: in Section \ref{section.1rob}, we recall the 2BSDEs framework and some useful results. Inspired by \cite{ekr} and \cite{him}, in Sections \ref{3}, \ref{fr}, \ref{4} and \ref{5} we solve the problem for robust exponential utility, robust power utility and robust logarithmic utility, which, unlike in \cite{dk}, are not bounded. Finally, in Section \ref{6}, we give some examples where we can explicitly solve the robust utility maximization problems by finding the solution of the associated 2BSDEs, and we provide some insights and comparisons with the classical dynamic programming approach adopted in the seminal work of Merton \cite{mer69}.

\vspace{0.8em}
\section{Preliminaries}
\label{section.1rob}
We will start by recalling some notations and notions related to the theory of 2BSDEs, which are the main tool in our approach to the robust utility maximization problem.

\subsection{Probability spaces}
Let $\Omega:=\left\{\omega\in C([0,T])^{d},\ \omega (0)=0\right\}$ be the canonical space, $B$ the canonical process, and $\mathbb{F}=\left(\mathcal F_t\right)_{0\leq t\leq T}$ the filtration generated by $B$. We will also make use of the right-limit $\mathbb F^+=\left(\mathcal F_{t^+}\right)_{0\leq t\leq T}$ of $\mathbb F$. Let $\mathbb{P}_{0}$ be the Wiener measure. As recalled in \cite{stz}, we can construct the quadratic variation of $B$ and its density $\widehat a$ pathwise.

\vspace{0.5em}
Let $\overline{\mathcal P}_W$ denote the set of all local martingale measures $\mathbb P$ such that $\mathbb P-a.s.$ for $t\in[0,T]$
\begin{equation}
\left<B\right>_t \text{ is absolutely continuous with respect to $t$ and $\widehat a$ takes values in $\mathbb S_d^{>0}$,}
\label{eq:}
\end{equation}
where $\mathbb S_d^{>0}$ denotes the space of all $d\times d$ real valued positive definite matrices. As in \cite{stz}, we concentrate on the subclass $\overline{\mathcal P}_S\subset\overline{\mathcal P}_W$ consisting of all probability measures

\begin{equation}
\mathbb P^\alpha:=\mathbb P_0\circ (X^\alpha)^{-1} \text{ where }X_t^\alpha:=\int_0^t\alpha_s^{1/2}dB_s,\text{ }t\in [0,T],\text{ }\mathbb P_0-a.s.,
\end{equation}
for some $\mathbb F$-progressively measurable process $\alpha$ in $\mathbb S_d^{>0}$ with $\int_0^T|\alpha_t|dt<+\infty$, $\mathbb {P}_0-a.s.$ Notice that $\alpha_s^{1/2}$ is just the square-root of the positive definite matrix $\alpha_s$.

\vspace{0.8em}
Finally, we fix $\underline{a}, \overline a \in\mathbb S_d^{>0}$ such that $\underline a\leq\overline a$ (for the usual order on positive definite matrices, i.e. $(\overline{a}-\underline{a})\in\mathbb S_d^{>0}$)  and we define the class
$$\mathcal P_H:=\left\{\mathbb P\in\overline{\mathcal P}_S \text{ s.t. } \underline a\leq \widehat a_t\leq \overline a, dt\times \mathbb P-a.e.\right\}$$

\vspace{0.8em}
which is a particular case of Definition $2.6$ in \cite{stz}, the main differences in our case being that the two bounds on $\widehat a$ are independent of the probability measures and that $\widehat F^0$ (introduced and defined below in Section \ref{pre}) is bounded. Throughout the paper it is assumed that $\mathcal{P}_{H}$ is not empty.

\vspace{0.8em}
For every $(t,\mathbb P)\in[0,T]\times\mathcal P_H$, we also define the class of probability measures which coincide with $\mathbb P$ up to $t^+$
$$\mathcal P_H(t^+,\mathbb P):=\left\{\mathbb P^{'}\in\mathcal P_H,\ \mathbb P^{'}=\mathbb P\text{ on }\mathcal F_{t^+}\right\}.$$

\vspace{0.8em}
\begin{Definition}
We say that a property holds $\mathcal{P}_{H}$-quasi-surely ($\mathcal{P}_{H}$-q.s. for short) if it holds $\mathbb{P}$-a.s. for all $\mathbb{P}\in\mathcal{P}_{H}$.
\end{Definition}

\begin{Remark}
The filtration $\mathbb F^+$ is right-continuous but not complete under each $\mathbb P\in\mathcal P_H$. Moreover, it is not possible to complete the filtration for each $\mathbb P$ since the measures are singular. It is of course a major drawback since many results of the general theory of processes rely on the fact that the underlying filtrations satisfy the usual hypotheses of right-continuity and completeness. However, this problem was solved in Lemma $2.4$ of \cite{stz3}, which implies that for every $\mathbb P\in\mathcal P_H$, we can always consider a version of our processes which is progressively measurable for the completion of $\mathbb F^+$ under $\mathbb P$.
\end{Remark}

\subsection{Spaces and Norms}\label{space}
We now recall from Possama\"i and Zhou \cite{posmz1} the spaces and norms which will be needed for the formulation of the quadratic second order BSDEs.

\vspace{0.8em}
$\mathbb L^{\infty}_H$ is the space of random variables which are bounded quasi-surely endowed with the norm
$$\No{\xi}_{\mathbb L^{\infty}_H}:=\underset{\mathbb{P} \in \mathcal{P}_H}{\sup}\No{\xi}_{L^{\infty}(\mathbb P)}.$$

\vspace{0.8em}
For $p\geq 1$, $\mathbb H^{p}_H$ denotes the space of all $\mathbb F^+$-progressively measurable $\mathbb R^d$-valued processes $Z$ with
$$\No{Z}_{\mathbb H^{p}_H}^p:=\underset{\mathbb{P} \in \mathcal{P}_H}{\sup}\mathbb E^{\mathbb P}\left[\left(\int_0^T|\widehat a_t^{1/2}Z_t|^2dt\right)^{\frac p2}\right]<+\infty.$$

\vspace{0.8em}
$\mathbb B\rm{MO}(\mathcal P_H)$ denotes the space of all $\mathbb F^+$-progressively measurable $\mathbb R^d$-valued processes $Z$ with
$$\No{Z}_{\mathbb B\rm{MO}(\mathcal P_H)}:=\underset{\mathbb{P} \in \mathcal{P}_H}{\sup}\No{\int_0^.Z_sdB_s}_{\rm{BMO}(\mathbb P)},$$
where $\No{\cdot}_{\rm{BMO}(\mathbb P)}$ is the usual $\rm{BMO}(\mathbb P)$ norm under $\mathbb P$, that is to say
$$\No{\int_0^.Z_sdB_s}_{\rm{BMO}(\mathbb P)}^2:=\underset{\tau\in\mathcal T_{0,T}}{\esup^\mathbb P}\No{\mathbb E^\mathbb P_\tau\left[\int_\tau^T\abs{\widehat a_tZ_t}^2dt\right]}_{L^\infty(\mathbb P)},$$
where $\mathcal T_{0,T}$ denotes the stopping times with value in $[0,T]$. We abuse notation and say that $\int_0^.Z_sdB_s$ is a $\mathbb B\rm{MO}(\mathcal P_H)$ martingale if $Z\in\mathbb B\rm{MO}(\mathcal P_H)$.

\vspace{0.5em}
$\mathbb D^{\infty}_H$ denotes the space of all $\mathbb F^+$-progressively measurable $\mathbb R$-valued processes $Y$ with
$$\mathcal P_H-q.s. \text{ c\`adl\`ag paths, and }\No{Y}_{\mathbb D^{\infty}_H}:=\underset{0\leq t\leq T}{\sup}\No{Y_t}_{L^{\infty}_H}<+\infty.$$

\vspace{0.5em}
We also recall the following space which is important in the formulation of Lipschitz 2BSDEs in \cite{stz}. For any $\kappa\in(1,2]$, $\mathbb L^{2,\kappa}_H$ is the space of random variables $\xi$ such that
$$\No{\xi}^2_{\mathbb L^{2,\kappa}_H}:= \underset{\mathbb{P} \in \mathcal{P}_H}{\sup}\mathbb E^\mathbb P\left[\underset{0\leq t\leq T}{\esup}^\mathbb P\ \underset{\mathbb P^{'}\in\mathcal P_H(t^+,\mathbb P)}{\esup^\mathbb P}\ \left(\mathbb E_t^{\mathbb P^{'}}\left[\abs{\xi}^\kappa\right]\right)^{\frac2\kappa}\right]<+\infty.$$

\vspace{0.5em}
Finally, we denote by $\mbox{UC}_b(\Omega)$ the collection of all bounded and uniformly continuous maps $\xi:\Omega\rightarrow \mathbb R$ with respect to the $\No{\cdot}_{\infty}$-norm, and we let
$$\mathcal L^{\infty}_H:=\text{the closure of $\mbox{UC}_b(\Omega)$ under the norm $\No{\cdot}_{\mathbb L^{\infty}_H}$},$$

and
$$\mathcal L^{2,\kappa}_H:=\text{the closure of $\mbox{UC}_b(\Omega)$ under the norm $\No{\cdot}_{\mathbb L^{2,\kappa}_H}$}.$$

\subsection{The quadratic generator}\label{pre}

We consider a map $H_t(\omega,z,\gamma):[0,T]\times\Omega\times\mathbb{R}^d\times D_H\rightarrow \mathbb{R}$, where $D_H \subset \mathbb{R}^{d\times d}$ is a given subset containing $0$.

\vspace{0.5em}
Define the corresponding conjugate of $H$ w.r.t. $\gamma$ by
\begin{align*}
&F_t(\omega,z,a):=\underset{\gamma \in D_H}{\rm{sup}}\left\{\frac12{\rm{Tr}} (a\gamma)-H_t(\omega,z,\gamma)\right\} \text{ for } a \in \mathbb S_d^{>0}\\[0.3em]
&\widehat{F}_t(z):=F_t(\omega, z,\widehat{a}_t) \text{ and } \widehat{F}_t^0:=\widehat{F}_t(0).
\end{align*}

\vspace{0.5em}
We denote by $D_{F_t(z)}$ the domain of $F$ in $a$ for a fixed $(t,\omega,z)$. As in \cite{posmz1}, the generator $F$ is supposed to verify either

\begin{Assumption} \label{assump.hquad}
\begin{itemize}
\item[\rm{(i)}] The domain $D_{F_t(z)}=D_{F_t}$ is independent of $(\omega,z)$.
\item[\rm{(ii)}] For fixed $(z,a)$, $F$ is $\mathbb{F}$-progressively measurable.
\item[\rm{(iii)}] $F$ is uniformly continuous in $\omega$ for the $||\cdot||_\infty$ norm.
\item[\rm{(iv)}] $F$ is continuous in $z$ and has the following growth property. There exists $(\alpha,\gamma)\in \mathbb R_+\times \mathbb R_+/\left\{0\right\}$ such that
$$\abs{ F_t(\omega,z,a)}\leq \alpha+\frac\gamma2\abs{ a^{1/2}z}^2,\text{ for all }(t,z,\omega,a).$$
\item[\rm{(v)}] $F$ is $C^2$ in $z$, and there are constants $r$ and $\theta$ such that for all $(t,\omega,z,a)$,
$$ \lvert D_z F_t(\omega,z,a)\rvert\leq r + \theta\abs{a^{1/2}z},\ \lvert D^2_{zz} F_t(\omega,y,z,a)\rvert\leq \theta.$$
\end{itemize}
\end{Assumption}
or
\begin{Assumption} \label{assump.hh}
Let points $\rm{(i)}$ through $\rm{(iv)}$ of Assumption \ref{assump.hquad} hold, and
\begin{itemize}
 \item[\rm{(v)}] $\exists \mu>0$ and a progressively measurable process $\phi\in\mathbb B\rm{MO}(\mathcal P_H)$ such that for all $(t,z,z',\omega,a),$
$$\abs{ F_t(\omega,z,a)- F_t(\omega,z',a)-\phi_t. a^{1/2}(z-z')}\leq \mu a^{1/2}\abs{z-z'}\left(|a^{1/2}z|+|a^{1/2}z'|\right).$$
\end{itemize}
\end{Assumption}

\vspace{0.8em}
\begin{Remark}
Notice that Assumption \ref{assump.hquad}(iv) implies that $\underset{0\leq t\leq T}{\esup}\abs{\widehat F^0_t}\in \mathbb L^\infty_H$.
\end{Remark}

\subsection{Quadratic 2BSDE}

In the sequel we will have to deal with the following type of 2BSDEs

\begin{equation}
\label{STZeq3.1}
Y_{t}=\xi-\int^{T}_{t}\widehat{F}_s(Z_{s})ds-\int^{T}_{t}Z_{s}dB_{s}+K_{T}-K_{t},\;0\leq t\leq T, \;\mathcal{P}_{H}-q.s.
\end{equation}
\begin{Definition}
Given $\xi\in\mathcal{L}^{\infty}_{H}$, we say $(Y,Z)\in\mathbb{D}^{\infty}_{H}\times\mathbb{H}^{2}_{H}$ is a solution to the 2BSDE \reff{STZeq3.1} if
\begin{itemize}
	\item $Y_{T}=\xi,\;\mathcal{P}_{H}-q.s.$
	
	\item For each $\mathbb{P}\in\mathcal{P}_{H}$, the process $K^{\mathbb{P}}$ defined below has nondecreasing paths, $\mathbb{P}$-a.s.
	\begin{equation}
	\label{STZeq3.2}
	K^{\mathbb{P}}_{t}:=Y_{0}-Y_{t}+\int^{t}_{0}\widehat{F}_{s}(Z_{s})ds+\int^{t}_{0}Z_{s}dB_{s},\;0\leq t\leq T,\;\mathbb{P}-a.s.
	\end{equation}

	\item The family of processes $\left\{K^{\mathbb{P}},\mathbb{P}\in\mathcal{P}_{H}\right\}$ defined in \reff{STZeq3.2} satisfies the following minimum condition
	\begin{equation}
	\label{STZeq3.3}
	K^{\mathbb{P}}_{t}=\underset{\mathbb{P}'\in\mathcal{P}_{H}(t^+,\mathbb{P})}{\einf^{\mathbb{P}}}E^{\mathbb{P}'}_{t}[K^{\mathbb{P'}}_{T}],\ \mathbb{P}-a.s.\text{ for all }\mathbb{P}\in\mathcal{P}_{H}, t\in[0,T].
	\end{equation}
%where for every $(t,\mathbb P)\in[0,T]\times\mathcal P_H$
%$$\mathcal P_H(t^+,\mathbb P):=\left\{\mathbb P^{'}\in\mathcal P_H,\ \mathbb P^{'}=\mathbb P\text{ on }\mathcal F_{t^+}\right\}.$$
\end{itemize}

\vspace{0.8em}
Moreover, if the family $\left\{K^{\mathbb{P}},\mathbb{P}\in\mathcal{P}_{H}\right\}$ can be aggregated into a universal process $K$, that is to say that for all $t\in[0,T]$
$$K_t=K_t^\mathbb P,\ \mathbb P-a.s.,\ \forall\mathbb P\in\mathcal P_H,$$
then we call $(Y,Z,K)$ a solution of the 2BSDE \reff{STZeq3.1}.
\end{Definition}

\vspace{0.8em}
Here one of the results proved in \cite{posmz1} is recalled (see Theorems $3.1$ and $4.1$)

\begin{Theorem}\label{exist}
Let $\xi\in\mathcal L^{\infty}_H$. Under Assumption \ref{assump.hquad} or Assumption \ref{assump.hh} with the addition that the norm of $\xi$ and the $\mathbb L^\infty_H$-norm of $\underset{0\leq t\leq T}{\esup}\ |\widehat F^0_t|$ are small enough, there is a unique solution $(Y,Z)\in\mathbb D^{\infty}_H\times\mathbb H^{2}_H$ of the $2BSDE$ \reff{STZeq3.1}. Moreover we have for all $t\in[0,T]$ and every $\mathbb P\in\mathcal P_H$
\begin{equation}\label{rep}
Y_t=\underset{\mathbb P^{'}\in\mathcal P_H(t^+,\mathbb P)}{\esup^\mathbb P} y_t^{\mathbb P^{'}},\ \mathbb P-a.s.,
\end{equation}
where $y^{\mathbb P^{'}}$ is the solution under $\mathbb P^{'}$ of the BSDE with generator $\widehat F$ and terminal condition $\xi$.

\end{Theorem}

\begin{Remark}
Assumption \ref{assump.hh} is weaker than Assumption \ref{assump.hquad}, but is sufficient to have existence of the quadratic $2$BSDE defined above only if the norms of the terminal condition $\xi$ and $\widehat F^0$ are small enough. Notice that since, for power and logarithmic utilities, the terminal condition will be equal to $0$, we only have a restriction on the norm of $\widehat F^0$ in these cases. We also emphasize that these restrictions are in no way necessary to obtain existence, but are artifacts of the type of proofs used in \cite{posmz1} to obtain existence of a 2BSDE with quadratic growth. Indeed, the proof relies at some point on the fact that solutions of standard BSDEs can be obtained through Picard iterations. Notice that this property was already needed in \cite{stz}. However, with a generator of quadratic growth, such a property was shown by Tevzadze \cite{tev} only if Assumption \ref{assump.hquad}$\rm{(v)}$ holds (see Proposition $2$), or if the terminal condition and $\widehat F^0$ are small enough and if Assumption \ref{assump.hh}$\rm{(v)}$ holds (see Proposition $1$). We conjecture that existence of solutions of 2BSDEs with quadratic growth should hold under less restrictive assumptions similar to those in \cite{elkarbar} (for instance $\xi$ would not need to be bounded and the generator would only need to be of quadratic growth), but this is left for future research.
\end{Remark}

\begin{Remark}
The representation \reff{rep} gives some insight into 2BSDEs. Since $Y$ can be written as a supremum of solution of BSDEs, we can interpret the increasing processes $K^\mathbb P$ as the instruments allowing $Y$ to remain above the corresponding $y^\mathbb P$. It is similar to reflected BSDEs with a lower obstacle. Moreover, the minimum condition \reff{STZeq3.3} tells us that this is done in a minimal way, making it the counterpart of the Skorokhod condition in our context.
\end{Remark}

\section{Robust utility maximization} \label{3}
We will now present the main problem of the paper and introduce a financial market with volatility uncertainty. The financial market consists of one bond with zero interest rate and $d$ stocks. The price process is given by
$$dS_{t}={\rm diag}\left[S_{t}\right](b_{t}dt+dB_{t}),\ \mathcal{P}_{H}-q.s.,$$
where $b$ is an $\R^d$-valued uniformly bounded stochastic process which is uniformly continuous in $\omega$ for the $||\cdot||_\infty$ norm.

\vspace{0.8em}
\begin{Remark}
The volatility is implicitly embedded in the model. Indeed, under each $\Pbb\in\Pset$, we have $dB_s\equiv\widehat{a}^{1/2}_tdW^{\Pbb}_t$ where $W^{\Pbb}$ is a Brownian motion under $\Pbb$. Therefore, $\widehat{a}^{1/2}$ plays the role of volatility under each $\Pbb$ and thus  makes it possible to model the volatility uncertainty. We also note that we make the uniform continuity assumption for $b$ to ensure that the generators of the $2$BSDEs obtained later satisfy Assumptions \ref{assump.hquad} or \ref{assump.hh}.
\end{Remark}

\vspace{0.8em}
We then denote $\pi=(\pi_t)_{0\leq t\leq T}$ a trading strategy, which is a $d$-dimensional $\F^+$-progressively measurable process, supposed to take its value in some closed set $A$. We refer to Definitions \ref{defexp}, \ref{defpow} and \ref{admiss.log} in the following sections for precise definitions of the set of admissible strategies $\mathcal A$ for the three utility functions studied. %We also denote $X^{\pi}_{t}$ the liquidation value at time $t$ of a trading strategy $\pi$ with positive initial capital $x$. Notice that we drop the $x$ in the notation of $X^{\pi}_{t}$ for simplification.

\vspace{0.8em}
The process $\pi^i_t $ describes the amount of money invested in stock $i$ at time $t$, with $1\leq i\leq d$. The number of shares is $\frac{\pi^i_t}{S^i_t} $. So the liquidation value of a trading strategy $\pi$ with positive initial capital $x$ is given by the following wealth process
$$X^{\pi}_{t}=x+\int^{t}_{0}\pi_{s}(dB_{s}+b_{s}ds),\ 0\leq t\leq T,\ \mathcal{P}_{H}-q.s.$$

%\begin{Remark}
Since zero interest rate was assumed, the amount of money in the bank $\pi^0$ does not appear in the wealth process $X$.
%\end{Remark}

\vspace{0.8em}
The problem of the investor in this financial market is to maximize the expected utility under model uncertainty of his terminal wealth $X^{\pi}_{T}-\xi$, where $\xi$ is a liability, that is to say a $\mathcal F_T$-measurable random variable. This liability could represent the value of any option or contract maturing at time $T$. It will always  be assumed that $\xi\in\mathcal L^{\infty}_H$.

\vspace{0.8em}
Denote by $U$ the utility function of the investor. The value function $V$ of the maximization problem therefore becomes
\begin{equation}
\label{optimprob0}	 V^{\xi}(x):= \underset{\pi\in\mathcal{A}}{{\rm sup}}\text{ }\underset{{\Q\in\Pset}}{{\rm inf}}\Ebb^{\Q}[U(X^{\pi}_T-\xi)].
\end{equation}

In the case where $\Pset$ contains only one probability measure, the problem reduces to the classical utility maximization problem.

\begin{Remark}
Due to the construction of 2BSDEs, we must have $\xi\in\mathcal L^{\infty}_H$. It is easy to see that $\xi$ can be constant, deterministic or in the form of $g(B_T)$ where $g$ is a Lipschitz bounded function, such as a Put or a Call spread payoff function. However, it can be noted that vanilla options payoffs with underlying $S$ may not be in $\mathcal L^{\infty}_H$. Indeed, we have in the one-dimensional framework
$$S_T=S_0\exp\left(\int^T_0 b_tdt-\frac{1}{2}\left\langle B\right\rangle_T + B_T\right),\ \mathcal{P}_{H}-q.s.$$

Since the quadratic variation of the canonical process can be written as follows
$$\underset{n\rightarrow+\infty}{\overline{\lim}}\sum_{i\leq 2^nt}\left(B_{\frac{i+1}{2^n}}(\omega)-B_{\frac{i}{2^n}}(\omega)\right)^2,$$
it is not too difficult to see that $S$ can be approximated by a sequence of random variables in $\rm{UC}_b(\Omega)$. Besides, this sequence converges in $\mathcal L^2_H$. However, we cannot be sure that it also converges in $\mathcal L^\infty_H$, which is the space of interest here.

\vspace{0.5em}
Of course, in the uncertain volatility framework, it seems to be a major drawback. Nevertheless, to deal with these options, it is sufficient to redo the whole $2$BSDE construction from scratch but taking the exponential of the Brownian motion under the Wiener measure as the canonical process instead of the Brownian motion itself. It would amount to restrict ourselves to the subset $\mathcal P_H^+$ of $\mathcal P_H$, containing only those $\mathbb P\in\mathcal P_H$ such that the canonical process is a positive continuous local martingale under $\mathbb P$.
\end{Remark}

To find the value function $V^{\xi}$ and an optimal trading strategy $\pi^*$, we follow the ideas of the general \textit{martingale optimality principle} approach as in \cite{ekr} and \cite{him}, but adapt it here to a nonlinear framework. Note that $\mathcal A$ is the admissibility set of the strategies $\pi$.

\vspace{0.5em}
Let $\{R^{\pi}\}_{\pi\in\mathcal A}$ be a family of processes which satisfy the following properties
\vspace{0.8em}
\begin{Properties}
\label{H0}
\begin{itemize}
	\item[\rm{(i)}] $R^{\pi}_{T}=U(X^{\pi}_{T}-\xi)$ for all $\pi\in\mathcal{A}$.
	
	\item[\rm{(ii)}] $R^{\pi}_{0}=R_{0}$ is constant for all $\pi\in\mathcal{A}$.
	
	\item[\rm{(iii)}] We have
	\begin{align*}
	 &R^{\pi}_{t}\geq\underset{\Pbb'\in\mathcal{P}_{H}(t^+,\Pbb)}{\einf^{\Pbb}}\Ebb^{\Pbb'}_t[R_T^\pi] , \text{ }\forall\pi\in\mathcal{A}\\[0.8em]
	&R^{\pi^*}_{t}=\underset{\Pbb'\in\mathcal{P}_{H}(t^+,\Pbb)}{\einf^{\Pbb}}\Ebb^{\Pbb'}_t[R_T^{\pi^*}] \text{ for some } \pi^*\in\mathcal{A}, \ \Pbb-a.s. \text{ for all } \Pbb\in\mathcal{P}_H.
	\end{align*}
\end{itemize}
\end{Properties}

\vspace{0.8em}
Then it follows that
\begin{equation}
\underset{\Pbb\in\mathcal{P}_H}{\rm inf}\Ebb^{\Pbb}[U(X^{\pi}_{T}-\xi)]\leq R_{0}=\underset{\Pbb\in\mathcal{P}_H}{\rm inf}\Ebb^{\Pbb}[U(X^{\pi^{*}}_{T}-\xi)]=V^{\xi}(x).
\end{equation}

\vspace{0.8em}
In the following sections we will follow the ideas of Hu, Imkeller and M\"uller \cite{him} to construct such a family for our three utility functions $U$.

\section{Robust exponential utility}\label{fr}
In this section,  the exponential utility function which is defined as
\begin{equation*}
U(x)=-\exp(-\beta x),\; x\in \mathbb{R} \text{ for $\beta>0$},
\end{equation*}
will be considered. In this context, the set of admissible trading strategies is defined as follows:

\begin{Definition}\label{defexp}

Let $A$ be a closed set in $\mathbb{R}^{d}$. The set of admissible trading strategies $\mathcal{A}$ consists of all $d$-dimensional progressively measurable processes, $\pi=(\pi_{t})_{0\leq t\leq T}$ satisfying
$$\pi\in\mathbb B\text{$\rm{MO}$}(\mathcal P_H) \text{ and }\pi_{t}\in A,\;dt\otimes \Pset-a.e.$$

\end{Definition}

\begin{Remark}\label{bilou}
Many authors have shed light on the natural link between BMO class, exponential uniformly integrable class and BSDEs with quadratic growth. See \cite{elkarbar}, \cite{bh} and \cite{him} among others. In the standard utility maximization problem studied in \cite{him}, their trading strategies satisfy a uniform integrability assumption on the family $\left(\exp(X^{\pi}_\tau)\right)_\tau$. Since the optimal strategy is a BMO martingale, it is easy to see that the utility maximization problem can also be solved if the uniform integrability assumption is replaced by a BMO assumption. However, at the end of the day, those two assumptions are deeply linked, as shown in the context of quadratic semimartingales in \cite{elkarbar}. Nonetheless, in our framework, as explained below in Remark \ref{bmoK}, it is necessary to generalize the BMO martingale assumption instead of the uniform integrability assumption. Moreover, as recalled in the Introduction, from a financial point of view these admissibility sets are related to absence of arbitrage in the market considered.
\end{Remark}

\subsection{Characterization of the value function and existence of optimal strategies}

\vspace{0.8em}
The investor wants to solve the maximization problem
\begin{eqnarray}
\label{valuefunction}	 V^{\xi}(x)\;:=\; \underset{\pi\in\mathcal{A}}{{\sup}}\text{ }\underset{{\Q\in\Pset}}{{\inf}}\Ebb^{\Q}\left[-\exp\left(X_T^\pi-\xi\right)\right].
\end{eqnarray}

\vspace{0.8em}
In order to construct a process $R^{\pi}$ which satisfies the Properties \ref{H0}, we set
\begin{equation*}
R^{\pi}_{t}=-\exp(-\beta(X^{\pi}_{t}-Y_{t})),\ t\in[0,T],\ \pi\in\mathcal{A},
\end{equation*}
where $(Y,Z)\in \mathbb D^{\infty}_H\times\mathbb H^{2}_H$ is the unique solution of a 2BSDE with a generator $\widehat{F}$ to be determined
\begin{equation*}
Y_{t}=\xi-\int^{T}_{t}Z_{s}dB_{s}-\int^{T}_{t}\widehat{F}(s,Z_{s})ds+K^{\Pbb}_{T}-K^{\Pbb}_{t},\ \Pbb-a.s.,\ \forall \Pbb\in\Pset.
\end{equation*}
%where the generator $\widehat{F}$ is chosen so that $R^{\pi}$ satisfies the Properties \ref{H0}.

\begin{Remark}
From Theorem $3.1$ of \cite{posmz1}, we have the following representation
$$Y_t=\underset{\Pbb'\in\mathcal{P}_{H}(t^+,\Pbb)}{\esup^{\Pbb}}y^{\Pbb'}_{t}.$$
Therefore, in general $Y_0$ is only $\mathcal F_{0^+}$-measurable and therefore not a constant. But by Proposition $4.2$ of \cite{posmz1}, we know that the process $Y$ is actually $\mathbb F$-measurable (it is true when the terminal condition is in $\rm{UC_b}(\Omega)$ and by passing to the limit when the terminal condition is in $\mathcal L^\infty_H$). This and the above representation easily imply that
$$Y_0=\underset{\Pbb'\in\mathcal{P}_{H}(0^+,\Pbb)}{\esup^{\Pbb}}y^{\Pbb'}_{0}=\underset{\Pbb'\in\mathcal{P}_{H}}{\sup}y^{\Pbb'}_{0}.$$
The Blumenthal Zero-One law then ensures that $Y_0$ is a constant.
\end{Remark}

Let us now define for all $a\in\mathbb S_d^{>0}$ such that $\underline a\leq a\leq \overline a$ the set $A_a$ by
$$A_a:=a^{1/2} A=\left\{a^{1/2}b,\ b\in A\right\}.$$

For any $a\in[\underline a,\overline a]$, the set $A_a$ is still closed. Moreover, since $ A\neq \varnothing$ we have
\begin{equation}\label{machin}
\min\left\{\abs{r},\text{ }r\in A_a\right\}\leq k,
\end{equation}
for some constant $k$ independent of $a$. We can now state the main result of this section

\begin{Theorem}\label{theoremexp}
Assume that $\xi\in\mathcal L^\infty_H$ and either that $\No{\xi}_{\mathbb L^{\infty}_H} + \underset{0\leq t\leq T}{\sup}\No{b_t}_{\mathbb L^\infty_H}$ is small and that $0\in A$, or that the set $A$ is $C^2$ (in the sense that its border is a $C^2$ Jordan arc). Then, the value function of the optimization problem \reff{valuefunction} is given by
\begin{equation*}
V^{\xi}(x)=-{{\rm exp}}\left(-\beta\left(x-Y_0\right)\right),
\end{equation*}
where $Y_0$ is defined as the initial value of the unique solution $(Y,Z)\in \mathbb D^{\infty}_H\times\mathbb H^{2}_H $ to the following $2$BSDE
\begin{equation}
\label{2BSDEQ}
Y_{t}=\xi-\int^{T}_{t}Z_{s}dB_{s}-\int^{T}_{t}\widehat{F}_s(Z_{s})ds+K^{\Pbb}_{T}-K^{\Pbb}_{t},\ \Pbb-a.s.,\ \forall \Pbb\in\Pset.
\end{equation}

\vspace{0.8em}
The generator is defined as follows
\begin{equation}
\label{generator1}
\widehat{F}_t(\omega,z):=F_t(\omega,z,\widehat{a}_t),
\end{equation}
where for all $t\in[0,T]$, $z\in\mathbb R^d$ and $a\in\mathbb S_d^{>0}$
$$F_t(\omega,z,a)=-\frac{\beta}{2}\text{$\rm{dist}$}^{2}\left(a^{1/2}z+\frac{1}{\beta}\theta_t(\omega),A_a\right)+z'a^{1/2}\theta_t(\omega)+\frac{1}{2\beta}\left|\theta_t(\omega)\right|^{2},$$
where $\theta_{t}(\omega):=a^{-\frac12}b_{t}(\omega)$ and for all $x\in\mathbb R^d$ and $E\subset\mathbb R^d$, $\text{$\rm{dist}$}(x,E)$ is the distance from $x$ to $E$.

\vspace{0.5em}
Moreover, there is an optimal trading strategy $\pi^*$ satisfying
\begin{equation}
\label{optim1} \widehat{a}_t^{1/2}\pi^*_t\in\Pi_{A_{\widehat a_t}}\left(\widehat{a}_t^{1/2}Z_t+\frac{1}{\beta}\widehat{\theta}_t\right), \ \ t\in[0,T],\text{ }\mathcal P_H-q.s.,
\end{equation}
with $\widehat \theta_t:=\widehat a^{-1/2}_tb_{t}$.
\end{Theorem}

\vspace{0.5em}
\proof
The proof is divided into $5$ steps. First, it is shown that the 2BSDE with the generator defined in (\ref{generator1}) has indeed a unique solution. Then, we prove a multiplicative decomposition for the process $R^\pi$ and some BMO integrability results on the process $Z$ and the optimal strategy $\pi^*$. Using these results, we are then able to show that $\rm{(iii)}$ of Properties \ref{H0} holds.

\vspace{0.8em}
\textbf{Step 1:} We show first that the 2BSDE \reff{2BSDEQ} has a unique solution. We need to verify that the generator $\widehat{F}$ satisfies the conditions of Assumption \ref{assump.hh} or \ref{assump.hquad}.

\vspace{0.5em}
First of all, $F$ defined above is a convex function of $a$, and thus for any $t\in[0,T]$, $F$ can be written as the Fenchel transform of a function
$$H_t(\omega,z,\gamma):=\underset{a \in D_F}{\sup}\left\{\frac12{\rm{Tr}}(a\gamma)-F_t(\omega,z,a)\right\} \text{ for } \gamma \in \mathbb{R}^{d\times d}.$$

\vspace{0.5em}
That $F$ satisfies the first two conditions of either Assumption \ref{assump.hh} or \ref{assump.hquad} is obvious. For Assumptions \ref{assump.hh}$\rm{(iii)}$ and \ref{assump.hquad}$\rm{(iii)}$, the assumption of boundedness and uniform continuity in $\omega$ on $b$ implies that $b^2$ is uniformly continuous in $\omega$. Since $b$ and $b^2$ are the only non-deterministic terms in $F$, then $F$ is also uniformly continuous in $\omega$.

\vspace{0.5em}
Then, since the distance function to a closed set is considered, we know that it is attained for some element of $\mathbb R^d$. Besides, as recalled earlier in \reff{machin}, there is a constant $k\geq 0$ such that
$$\min\left\{\left|d\right|,\ d\in A_{\widehat a_t}\right\}\leq k,\ dt\otimes \Pbb-a.e. \text{, for all }\mathbb P\in\mathcal P_H.$$

Then we get, for all $z\in\mathbb{R}^d,\;t\in[0,T] $,
$$\text{$\rm{dist}$}^2\left(\widehat{a}_t^{1/2}z+\frac{1}{\beta}\widehat{\theta}_t, A_{\widehat a_t}\right)\leq 2\left|\widehat{a}_t^{1/2}z\right|^2+2\left(\frac{1}{\beta}\left|\widehat{\theta}_t\right|+k\right)^2. $$

\vspace{0.5em}
Thus, we obtain from the boundedness of $\widehat{\theta} $
$$\left|\widehat{F}_t(z)\right|\leq c_0+c_1\left|\widehat{a}_t^{1/2}z\right|^2,$$
that is to say that Assumptions \ref{assump.hh}$\rm{(iv)}$ and \ref{assump.hquad}$\rm{(iv)}$ are satisfied.

\vspace{0.5em}
Finally, Assumption \ref{assump.hh}$\rm{(v)}$ is clear from the Lipschitz property of the distance function, and Assumption \ref{assump.hquad}$\rm{(v)}$ is also clear from the regularity assumption on the border of $A$.

\vspace{0.8em}
The terminal condition $\xi$ is in $\mathcal L^{\infty}_H$ and we have proved that the generator $\widehat{F}$ satisfies Assumption \ref{assump.hh} or Assumption \ref{assump.hquad}. Moreover, by definition of $F$, it is clear that if $b$ has a small $\mathbb L^\infty_H$-norm and if $0\in A$, then $\widehat F^0$ also has a small $\mathbb L^\infty_H$-norm. Indeed, we have
$$\widehat F_t^0=-\frac\beta2\text{$\rm{dist}$}\left(\frac{\theta_t}{\beta},A_{\widehat a_t}\right)+ \frac{1}{2\beta}\abs{\theta_t}^2,$$
which tends to $0$ as $b_t$ and thus $\theta_t$ goes to $0$ (it is clear for the second term on the right-hand side, and for the first, continuity of the distance function and the fact $0\in A$ ensure the result).

\vspace{0.8em}
Therefore Theorem \ref{exist} states that the $2$BSDE \reff{2BSDEQ} has a unique solution in $\mathbb D^{\infty}_H\times\mathbb H^{2}_H$.

\vspace{0.5em}
\textbf{Step 2:} We first decompose $R^{\pi}$ as the product of a process $M^{\pi}$ and a non-increasing process $N^{\pi}$ which is constant for some $\pi^{*}\in\mathcal{A}$.

\vspace{0.5em}
Define for all $\mathbb P\in\mathcal P_H$ any for any $t\in[0,T]$
\begin{equation*}
M^{\pi}_{t}=e^{-\beta(x-Y_{0})}\exp\left(-\int^{t}_{0}\beta(\pi_{s}-Z_{s})dB_{s}-\frac{1}{2}\int^{t}_{0}\beta^{2}\abs{\widehat{a}_s^{1/2}(\pi_{s}-Z_{s})}^{2}ds-\beta K^{\Pbb}_t\right),\ \Pbb-a.s.
\end{equation*}

We can then write for all $t\in[0,T]$
$$R^{\pi}_t= M^{\pi}_tN^{\pi}_t,$$
with $$N^{\pi}_{t}=-\exp\left(\int^{t}_{0}v(s,\pi_{s},Z_{s})ds\right),$$
and $$v(t,\pi,z)\;  \;=-\beta \pi b_{t}+\beta \widehat{F}_t(z)+\frac{1}{2}\beta^{2}\left|\widehat{a}_t^{1/2}\left(\pi-z\right)\right|^{2}.$$

Clearly, for every $t\in[0,T], $   $v(t,\pi_{t},Z_{t})$ can be rewritten in the following form
\begin{eqnarray*}
\frac{1}{\beta}v(t,\pi_{t},Z_{t})&=& \frac{\beta}{2}\left|\widehat{a}_t^{1/2}\pi_{t}\right|^{2}-\beta \pi^{'}_{t}\widehat{a}_t^{1/2}\left(\widehat{a}_t^{1/2}Z_{t}+\frac{1}{\beta}\widehat{\theta_{t}}\right)+\frac{\beta}{2}\left|\widehat{a}_t^{1/2}Z_{t}\right|^{2}+\widehat{F}_t(Z_{t})\\
&=&\frac{\beta}{2}\left|\widehat{a}_t^{1/2}\pi_{t}-\left(\widehat{a}_t^{1/2}Z_{t}+\frac{1}{\beta}\widehat{\theta_{t}}\right)\right|^{2}-Z^{'}_{t}\widehat{a}_t^{1/2}\widehat{\theta_{t}}-\frac{1}{2\beta}\left|\widehat{\theta_{t}}\right|^{2}+\widehat{F}_t(Z_{t}).
\end{eqnarray*}

By a classical measurable selection theorem (see \cite{elk} or Lemma $3.1$ in \cite{elkaroui}), we can define a progressively measurable process $\pi^*$ satisfying \reff{optim1}. Then, it follows from the definition of $\widehat{F}$ that $\mathcal P_H-q.s.$
\begin{itemize}
	\item $v(t,\pi_{t},Z_{t})\geq 0$ for all $\pi\in\mathcal{A}$, $t\in[0,t]$.
	\item $v(t,\pi^{*}_{t},Z_{t})=0$, $t\in[0,T]$,
\end{itemize}
which implies that the process $N^{\pi}$ is always non-increasing for all $\pi$ and is equal to $-1$ for $\pi^*$.

\vspace{0.8em}
\textbf{Step 3:}
In this step, we show that the processes
$$\int^{\cdot}_0Z_sdB_s,\ \ \int^{\cdot}_0\pi^*_sdB_s, $$
are $\mathbb B\rm{MO}(\mathcal P_H)$ martingales.

\vspace{0.8em}
First of all, by Lemma $2.1$ in \cite{posmz1}, we know that $\int^{\cdot}_0 Z_s dB_s $ is a $\mathbb B\rm{MO}(\mathcal P_H)$ martingale. By the triangle inequality and the definition of $\pi^*$ together with \reff{machin}, we have for all $t\in[0,T]$
\begin{eqnarray*}
\left|\widehat{a}_t^{1/2}\pi_t^*\right|& \leq & \left|\widehat{a}_t^{1/2}Z_t+\frac{1}{\beta}\widehat{\theta_t}\right| +\left|\widehat{a}_t^{1/2}\pi_t^*-\left(\widehat{a}_t^{1/2}Z_t+\frac{1}{\beta}\widehat{\theta}_t\right)\right|\\
&\leq & 2\left|\widehat{a}_t^{1/2}Z_t\right|+\frac{2}{\beta}\left|\widehat{\theta}_t\right|+k\leq 2\left|\widehat{a}_t^{1/2}Z_t\right|+k_1,
\end{eqnarray*}
where $k_1$ is a bound on $\widehat{\theta}$. Then, for every probability $\P\in\Pset$ and every stopping time $\tau\leq T $,
$$\Ebb_{\tau}^{\P}\left[\int^{T}_{\tau} \left|\widehat{a}_t^{1/2}\pi^*_t\right|^2 dt\right]\leq \Ebb_{\tau}^{\P}\left[\int^{T}_{\tau} 8\left|\widehat{a}_t^{1/2}Z_t\right|^2 dt+2Tk^2_1 \right],$$

\vspace{0.8em}
and therefore
%$$\No{\int^{\cdot}_0 \pi^*_s dB_s }_{\mathbb B\rm{MO}(\mathcal P_H)}\leq 8\No{\int^{\cdot}_0 Z_s dB_s }_{\mathbb B\rm{MO}(\mathcal P_H)}+2Tk^2_1.$$
$$\No{\pi^*}_{\mathbb B\rm{MO}(\mathcal P_H)}^2\leq 8\No{Z}_{\mathbb B\rm{MO}(\mathcal P_H)}^2+2Tk^2_1,$$
which implies the $\mathbb B\rm{MO}(\mathcal P_H)$ martingale property of $\int^{\cdot}_0 \pi^*_s dB_s $ as desired.

\vspace{1em}
\textbf{Step 4:} We then prove that $\pi^*\in\mathcal A$ and $R^{\pi^*}\equiv -M^{\pi^*}$ satisfies $\rm{(iii)}$ of Properties \ref{H0}, that is to say for all $t\in[0,T]$
$$\underset{\Pbb'\in\mathcal P_H(t^+,\mathbb P)}{\esup^{\Pbb}}\Ebb^{\Pbb'}_t\left[M^{\pi^*}_T\right]= M^{\pi^*}_t,\ \ \Pbb-a.s., \ \forall \Pbb\in\Pset.$$

\vspace{0.5em}
For a fixed $\Pbb'\in\mathcal P_H(t^+,\mathbb P)$, we denote $$L_t:=\int^{t}_{0}\beta(\pi^*_{s}-Z_{s})dB_{s}+\frac{1}{2}\int^{t}_{0}\beta^{2}\abs{\widehat{a}_s^{1/2}(\pi^*_{s}-Z_{s})}^{2}ds+\beta K^{\Pbb'}_t,\ 0\leq t\leq T,$$
then with It\^o's formula, we obtain for every $t\in[0,T]$, thanks to the $\mathbb B\rm{MO}(\mathcal P_H)$ property proved in Step $3$
\begin{align}\label{truc}
\nonumber\Ebb^{\Pbb'}_t\left[M^{\pi^*}_T\right]-M^{\pi^*}_t&=-\beta\Ebb^{\Pbb'}_t\left[\int^{T}_t M^{\pi^*}_{s^-} dK^{\Pbb'}_s\right]\\
&\hspace{0.9em}+\mathbb E_t^{\mathbb P^{'}}\left[\sum_{t\leq s\leq T}e^{- L_s}-e^{- L_{s^{-}}}+ e^{- L_{s^{-}}}(L_s-L_{s^{-}})\right].
\end{align}

First, we prove
$$\underset{\Pbb'\in\mathcal P_H(t^+,\mathbb P)}{\einf^{\Pbb}} \Ebb^{\Pbb'}_t\left[\int^{T}_t M^{\pi^*}_{s^-} dK^{\Pbb'}_s\right]=0,\ t\in[0,T],\ \Pbb-{\rm a.s.} $$

\vspace{0.8em}
For every $t$ and every $\Pbb'\in\mathcal P_H(t^+,\mathbb P)$, we have
$$0\leq\Ebb^{\Pbb'}_t\left[\int^{T}_t M^{\pi^*}_{s^-} dK^{\Pbb'}_s\right]\leq \Ebb^{\Pbb'}_t\left[\left(\underset{0\leq s\leq T}{{\rm sup}} M^{\pi^*}_s\right) \left(K^{\Pbb'}_T-K^{\Pbb'}_t\right)\right].$$
%$$0\leq\Ebb^{\Pbb'}_t\left[\int^{T}_t M^{\pi^*}_{s^-} dK^{\Pbb'}_s\right]\leq M^{\pi^*}_t\Ebb^{\Pbb'}_t\left[\left(\underset{t\leq s\leq T}{{\rm sup}} \left( M^{\pi^*}_t\right)^{-1} M^{\pi^*}_s\right) \left(K^{\Pbb'}_T-K^{\Pbb'}_t\right)\right].$$

\vspace{0.5em}
Besides, since $K^{\mathbb P^{'}}$ is non-decreasing, we obtain for all $s\geq t$
$$ M^{\pi^*}_s\leq e^{-\beta(x-Y_0)}\mathcal{E}\left(\beta\int_0^s \left(Z_u-\pi^*_u\right)dB_u \right).$$
%$$\left( M^{\pi^*}_t\right)^{-1} M^{\pi^*}_s\leq \mathcal{E}\left(\beta\int_t^s \left(Z_u-\pi^*_u\right)dB_u \right).$$

\vspace{1em}
Then, again thanks to Step $3$, we know that
%$$\int^{.}_0 \left(Z_s-\pi^*_s\right)dB_s \in \mathbb B\rm{MO}(\mathcal P_H),$$
$$ \left(Z_s-\pi^*_s\right) \in \mathbb B\rm{MO}(\mathcal P_H),$$
and thus the exponential martingale above is a uniformly integrable martingale for all $\mathbb P$ and is in $L^r_H$ for some $r>1$ (see Lemma $2.2$ in \cite{posmz1}). Thus, by H\"older inequality, we have for all $t\in[0,T]$
\begin{align*}
\Ebb^{\Pbb'}_t\left[\int^{T}_t M^{\pi^*}_{s^-} dK^{\Pbb'}_s\right]\leq e^{\beta(Y_0-x)} \Ebb^{\Pbb'}_t\left[\underset{0\leq s\leq T}{{\rm sup}} \mathcal{E}^r\left(\beta\int_0^s \left(Z_u-\pi^*_u\right)dB_u\right)\right]^{\frac1r}\Ebb^{\Pbb'}_t\left[\left(K^{\Pbb'}_T-K^{\Pbb'}_t\right)^q\right]^{\frac1q}.
\end{align*}
%\begin{align*}
%\Ebb^{\Pbb'}_t\left[\int^{T}_t M^{\pi^*}_{s^-} dK^{\Pbb'}_s\right]\leq M^{\pi^*}_t \Ebb^{\Pbb'}_t\left[\underset{t\leq s\leq T}{{\rm sup}} \mathcal{E}^r\left(\beta\int_t^s \left(Z_u-\pi^*_u\right)dB_u\right)\right]^{\frac1r}\Ebb^{\Pbb'}_t\left[\left(K^{\Pbb'}_T-K^{\Pbb'}_t\right)^q\right]^{\frac1q}.
%\end{align*}

\vspace{0.8em}
With Doob's maximal inequality, we have for every $t\in[0,T]$
\begin{eqnarray*}
\Ebb^{\Pbb'}_t\left[\underset{0\leq s\leq T}{{\rm sup}} \mathcal{E}^r\left(\beta\int_0^s \left(Z_u-\pi^*_u\right)dB_u\right)\right]^{1/r}&\leq& C\Ebb^{\Pbb'}_t\left[\mathcal{E}^r\left(\beta\int_0^T \left(Z_u-\pi^*_u\right)dB_u\right)\right]^{1/r}< +\infty,
\end{eqnarray*}
where $C$ is an universal constant which can change value from line to line.

\vspace{0.8em}
Then by the Cauchy-Schwarz inequality, we get for $0\leq t\leq T$
\begin{align*}
\Ebb^{\Pbb'}_t\left[\left(K^{\Pbb'}_T-K^{\Pbb'}_t\right)^q\right]^{1/q}&\leq C\left(\Ebb^{\Pbb'}_t\left[\left(K^{\Pbb'}_T-K^{\Pbb'}_t\right)\right]\Ebb^{\Pbb'}_t\left[\left(K^{\Pbb'}_T-K^{\Pbb'}_t\right)^{2q-1}\right]\right)^{\frac{1}{2q}}\\
&\leq C\left(\underset{\Pbb'\in\mathcal P_H(t^+,\mathbb P)}{\esup^{\Pbb}}\Ebb^{\Pbb'}_t\left[\left(K^{\Pbb'}_T-K^{\Pbb'}_t\right)^{2q-1}\right]\right)^{\frac{1}{2q}}\left(\Ebb^{\Pbb'}_t\left[\left(K^{\Pbb'}_T-K^{\Pbb'}_t\right)\right]\right)^{\frac{1}{2q}}.
\end{align*}

Arguing as in the proof of Theorem $3.1$ in \cite{posmz1} we know that
$$\left(\underset{\Pbb'\in\mathcal P_H(t^+,\mathbb P)}{\esup^{\Pbb}}\Ebb^{\Pbb'}_t\left[\left(K^{\Pbb'}_T-K^{\Pbb'}_t\right)^{2q-1}\right]\right)^{\frac{1}{2q}}<+\infty, \ 0\leq t\leq T.$$

Hence, we obtain for $0\leq t\leq T$
\begin{eqnarray*}
0\leq \underset{\Pbb'\in\mathcal P_H(t^+,\mathbb P)}{\einf^{\Pbb}} \Ebb^{\Pbb'}_t \left[\int^{T}_t M^{\pi^*}_{s^-} dK^{\Pbb'}_s\right] \leq C \underset{\Pbb'\in\mathcal P_H(t^+,\mathbb P)}{\einf^{\Pbb}} \left(\Ebb^{\Pbb'}_t\left[\left(K^{\Pbb'}_T-K^{\Pbb'}_t\right)\right]\right)^{\frac{1}{2q}}=0,
\end{eqnarray*}
which means
$$\underset{\Pbb'\in\mathcal P_H(t^+,\mathbb P)}{\einf^{\Pbb}}\Ebb^{\Pbb'}_t\left[\int^{T}_t M^{\pi^*}_{s^-} dK^{\Pbb'}_s\right]=0,\ 0\leq t\leq T,\ \mathbb P-a.s.$$

\vspace{0.8em}
Finally, we have for every $t\in[0,T]$
\begin{align*}
&\underset{\Pbb'\in\mathcal P_H(t^+,\mathbb P)}{\einf^{\Pbb}}\Ebb^{\Pbb'}_t\left[\int^{T}_t M^{\pi^*}_{s^-} dK^{\Pbb'}_s
-\sum_{t\leq s\leq T}\exp(-\beta L_s)-\exp(-\beta L_{s^{-}})+\beta \exp(-\beta L_{s^{-}})(L_s-L_{s^{-}})\right]\\
&\leq \underset{\Pbb'\in\mathcal P_H(t^+,\mathbb P)}{\einf^{\Pbb}}\Ebb^{\Pbb'}_t\left[\int^{T}_t M^{\pi^*}_{s^-} dK^{\Pbb'}_s\right]\\
&\hspace{0.9em}-\underset{\Pbb'\in\mathcal P_H(t^+,\mathbb P)}{\einf^{\Pbb}}\Ebb^{\Pbb'}_t\left[\sum_{t\leq s\leq T}\exp(-\beta L_s)-\exp(-\beta L_{s^{-}})+\beta \exp(-\beta L_{s^{-}})(L_s-L_{s^{-}})\right]\\&\leq 0,\ \mathbb P-a.s.,
\end{align*}
because the function $x\rightarrow \exp(- x)$ is convex and the jumps of $L$ are positive. Hence, using \reff{truc}, we have for every $t\in[0,T]$
$$\underset{\Pbb'\in\mathcal P_H(t^+,\mathbb P)}{\esup^{\Pbb}}\Ebb^{\Pbb'}_t\left[M_T^{\pi^*}-M_t^{\pi^*}\right]\geq 0,\ \mathbb P-a.s.$$

\vspace{0.5em}
But by definition $M^{\pi^*}$ is the product of a martingale and a positive non-increasing process and is therefore a supermartingale. It implies that for every $t\in[0,T]$
$$\underset{\Pbb'\in\mathcal P_H(t^+,\mathbb P)}{\esup^{\Pbb}}\Ebb^{\Pbb'}_t\left[M_T^{\pi^*}-M_t^{\pi^*}\right]= 0,\ \mathbb P-a.s.$$

Finally, $\pi^*$ is an admissible strategy, $R^{\pi^*} $ satisfies $\rm{(iii)}$ of Properties \ref{H0} and
\begin{eqnarray*}
R^{\pi^*}_0&=&\underset{\Pbb\in\Pset}{{\rm inf}}\Ebb^{\Pbb}\left[-\exp\left(-\beta\left(x+\int^{T}_{0}\pi^*_s\left(dB_s+\theta_s ds\right)-\xi\right)\right) \right]\\
&=&-\exp\left(-\beta\left(x-Y_0\right) \right).
\end{eqnarray*}

\vspace{1em}
\textbf{Step 5: } Next we will show that for all $\pi\in\mathcal{A} $, $R^{\pi}$ satisfies $\rm{(iii)}$ of Properties \ref{H0}, that is, for every $t\in[0,T]$
$$\underset{\Pbb'\in\mathcal{P}_H(t^+,\Pbb)}{\einf^{\Pbb}}\Ebb^{\Pbb'}_t[-\exp(-\beta(X^{\pi}_{T}-\xi))]\leq R^{\pi}_{t}, \ \Pbb-a.s.$$

\vspace{0.8em}
Since $\pi\in\mathcal{A}$, the process
$$\int^{.}_0 \left(Z_s-\pi_s\right)dB_s,$$
is a $\mathbb B\rm{MO}(\mathcal P_H)$ martingale. Then the process
$$G^{\pi}=\exp\left(-\beta(x-Y_0)\right)\mathcal{E}\left(-\beta\int_0^.\left(\pi_s-Z_s\right)dB_s \right),$$
is a uniformly integrable martingale under each $\Pbb\in\Pset$.

\vspace{0.8em}
As in the previous steps, we write $R^{\pi}$ as $R^{\pi}=M^{\pi}N^{\pi}$, where $N^{\pi}$ is a negative non-increasing process. We then have for $0\leq s\leq t\leq T$
\begin{align*}
\underset{\Pbb'\in\mathcal P_H(s+,\mathbb P)}{\einf^{\Pbb}}\Ebb_s^{\Pbb'}\left[M^{\pi}_{t}N^{\pi}_{t}\right] &\leq \underset{\Pbb'\in\mathcal P_H(s+,\mathbb P)}{\einf^{\Pbb}}\Ebb_s^{\Pbb'}\left[M^{\pi}_{t}N^{\pi}_{s}\right]\\
&=\underset{\Pbb'\in\mathcal P_H(s+,\mathbb P)}{\esup^{\Pbb}}\Ebb_s^{\Pbb'}\left[M^{\pi}_{t}\right]N^{\pi}_{s},\text{ }\mathbb P-a.s.
\end{align*}
because $N^{\pi}$ is negative. By the same arguments as in Step $3$ for $M^{\pi^{*}}$, we have for $0\leq s\leq t\leq T$
$$\underset{\Pbb'\in\mathcal P_H(s+,\mathbb P)}{\esup^{\Pbb}}\Ebb_s^{\Pbb'}\left[M^{\pi}_{t}\right]=M^{\pi}_{s},\ \Pbb-a.s.$$

\vspace{0.8em}
Therefore the following inequality holds for $0\leq s\leq t\leq T$
$$\underset{\Pbb'\in\mathcal P_H(s+,\mathbb P)}{\einf^{\Pbb}}\Ebb_s^{\Pbb'}\left[R^{\pi}_{t}\right] \leq R^{\pi}_{s},\ \Pbb-a.s. $$
which ends the proof.
\ep

\vspace{0.8em}
\begin{Remark}\label{bmoK}
Here, it can be seen why it is essential in this context to have strong integrability assumptions on the trading strategies. Indeed, in the proof of the above property for $M^{\pi^*}$, the fact that the stochastic integral
$$\int_0^\cdot\pi_s^*dB_s,$$
is a $\mathbb B\rm{MO}(\mathcal P_H)$ martingale allowed us to control the moments of its stochastic exponential, which in turn allowed us to deduce from the minimal property for $K^\mathbb P$ a similar minimal property for
$$\int_0^\cdot M_s^{\pi^*}dK_s^\mathbb P.$$

\vspace{0.8em}
This term is new when compared with the context of \cite{him}. To deal with it, the $\mathbb B\rm{MO}(\mathcal P_H)$ property has to be imposed. Note however that since the optimal strategy already has that property, we do not lose much by restricting the strategies.
\end{Remark}

\begin{Remark}
We note that the approach still works when there are no constraints on trading strategies. %as studied in \cite{pham}.
In this case, the 2BSDE related to the maximization problem has a uniformly Lipschitz generator, and we are in the context of complete markets. Then, the theory developed in \cite{stz} for Lipschitz 2BSDEs can also be used.
\end{Remark}

\subsection{A min-max property}

By comparing the value function of our robust utility maximization problem and the one presented in \cite{him} for standard utility maximization problem, we are able to obtain a min-max property similar to that obtained by Denis and Kervarec in \cite{dk}. We observe that we were only able to prove this property after having solved the initial problem, unlike in the approach of \cite{dk}.

\begin{Theorem}\label{minmax}
Under the previous assumptions on the probability measures set $\Pset$ and the admissible strategies set $\Ac$, the following min-max property holds.
$$\underset{\pi\in\mathcal{A}}{{\rm sup}}\text{  }\underset{\P\in\Pset}{{\rm inf}}\Ebb^{\P}\left[ R^{\pi}_T\right]= \underset{\P\in\Pset}{{\rm inf}}\text{ }\underset{\pi\in\mathcal{A}}{{\rm sup}}\text{ }\Ebb^{\P}\left[ R^{\pi}_T\right]= \underset{\P\in\Pset}{{\rm inf}}\text{ }\underset{\pi\in\mathcal{A}^{\P}}{{\rm sup}}\Ebb^{\P}\left[ R^{\pi}_T\right],  $$
where $\mathcal{A}^{\P}$ is the set consisting of trading strategies $\pi$ which are in $A$ and such that the process $\int_0^.\pi_sdB_s$ is a BMO($\mathbb P$) martingale.
\end{Theorem}

\proof
First note that we have $$D:=\underset{\pi\in\mathcal{A}}{{\rm sup}}\text{  }\underset{\P\in\Pset}{{\rm inf}}\Ebb^{\P}\left[ R^{\pi}_T\right]\leq \underset{\P\in\Pset}{{\rm inf}}\text{ }\underset{\pi\in\mathcal{A}}{{\rm sup}}\text{ }\Ebb^{\P}\left[ R^{\pi}_T\right]\leq \underset{\P\in\Pset}{{\rm inf}}\text{ }\underset{\pi\in\mathcal{A}^{\P}}{{\rm sup}}\Ebb^{\P}\left[ R^{\pi}_T\right]=:C.$$

Indeed, the first inequality is obvious and the second one follows from the fact that for all $\mathbb P$, $\mathcal{A}\subset\mathcal{A}^\mathbb P$.

\vspace{0.8em}
That $C\leq D$ remains to be proved. By the previous sections, we know that
%$$D=-{\rm exp}\left(-\beta\left(x-Y_0\right)\right).$$
\begin{equation*}
D=-{\rm exp}\left(-\beta\left(x-Y_0\right)\right).
\end{equation*}

Moreover, we know from \cite{posmz1} that we have a representation for $Y_0$,
\begin{equation*}
Y_0=\underset{\P\in\Pset}{{\rm sup}} y^{\P}_0,
\end{equation*}
where $y^{\P}_0$ is the solution of  the standard BSDE with the same generator $\widehat{F}$. On the other hand, it can be observed from \cite{him} that $$C= \underset{\P\in\Pset}{{\rm inf}} \left[-{\rm exp}\left(-\beta\left(x-y^{\P}_0\right)\right)\right],$$
implying that $C=D$.
\ep

%\subsection{Existence of an optimal probability measure}
%In general, we are not sure whether an optimal probability measure for the optimization problem \reff{valuefunction} exists or not. We will %present a case for which we were able to prove this existence. Throughout this subsection, we assume that we restrict ourselves to subset of $%\mathcal P_H$ which satisfies the conditions of Theorem $5.2$ of \cite{posmz1}.
%\vspace{0.8em}
%\begin{Proposition}\label{probopt}
%Assume that the liability $\xi$ is in $\rm{UC_b(\Omega)}$. Then an optimal probability measure $\Pbb^*$ exists.
%\end{Proposition}
%\proof
%We prove in the appendix that under this assumption, there exists a probability measure $\mathbb P^*$ such that
%$$Y_0= y_0^{\mathbb P^*}.$$
%\vspace{0.8em}
%From our min-max property, it is then clear that the measure $\mathbb P^*$ is optimal.
%\ep
%\vspace{0.8em}
\subsection{Indifference pricing via robust utility maximization} %in the Uncertain Volatility Model
It has been shown in \cite{ekr} that in a market model with constraints on the portfolios, if the indifference price for a claim $\Phi$  is defined as the smallest number $p$ such that
$$\underset{\pi}{\sup}\text{ } \Ebb\left[-\exp\left(-\beta \left(X^{x+p,\pi}-\Phi\right)\right)\right] \geq \underset{\pi}{\sup}\text{ } \Ebb\left[-\exp\left(-\beta X^{x,\pi}\right)\right],$$
where $X^{x,\pi}$ is the wealth associated with the portfolio $\pi$ and initial value $x$, then this problem turns into the resolution of a BSDE with quadratic growth generator.

\vspace{0.8em}
In this framework of uncertain volatility, the problem of indifference pricing of a contingent claim $\Phi$ boils down to solving the following equation in $p$
$$V^{0}(x)=V^{\Phi}(x+p).$$

Thanks to our results, we know that if $\Phi\in\mathcal L^\infty_H$ then the two sides of the above equality can be calculated by solving $2$BSDEs. Price $p$ can therefore be calculated as soon as the $2$BSDEs  can be solved (explicitly or numerically). Two examples are provided in Section \ref{6}.

\section{Robust power utility}\label{4}
In this section, we will consider the power utility function
$$U(x)=-\frac{1}{\gamma}x^{-\gamma},\ x > 0,\ \gamma > 0.$$%\in(0,1)

\vspace{0.8em}
Here  a different notion of trading strategy will be used: $\rho=(\rho^i)_{i=1,\ldots,d} $ denotes the proportion of wealth invested in stock $i$. The number of shares of stock $i$ is given by $\frac{\rho^i_tX_t}{S^i_t}$.

\vspace{0.8em}
Then the wealth process is defined as
\begin{equation}
X^{\rho}_t=x+\int^{t}_{0}{\sum^{d}_{i=1}}\frac{X^{\rho}_s {\rho}^i_{s}}{S^i_{s}} dS^i_{s}=x+\int^{t}_{0}X^{\rho}_{s}{\rho}_s\left(dB_s+b_sds\right),\ \mathcal{P}_{H}-q.s.
\end{equation}
and the initial capital $x$ is positive.

\vspace{0.8em}
In the present setting, the set of admissible strategies is defined as follows

\begin{Definition}\label{defpow}
Let $A$ be a closed set in $\mathbb{R}^{d}$. The set of admissible trading strategies $\mathcal{A}$ consists of all $\mathbb{R}^d$-valued progressively measurable processes $\rho=(\rho_{t})_{0\leq t\leq T}$ satisfying
$$\rho\in\mathbb B\text{$\rm{MO}$}(\mathcal P_H)\ \text{and}\ \rho\in A, \ dt\otimes \mathcal P_H-a.e.$$
\end{Definition}

\vspace{0.8em}
The wealth process $X^{\rho}$ can be written as
$$X^{\rho}_t=x\mathcal{E}\left(\int_0^t \rho_s(dB_s+b_sds)\right), \ t\in\left[0,T\right],\ \mathcal{P}_{H}-q.s.$$

\vspace{0.8em}
Then for every $\rho\in\mathcal{A} $, the wealth process $X^{\rho} $ is a local $\P $-martingale bounded from below, hence, a $\P $-supermartingale, for all $\P\in\Pset$.

\vspace{0.8em}
We suppose that there is no liability ($\xi=0$). Then the investor faces the maximization problem
\begin{equation}
\label{robustpower}
V(x)=\underset{\rho\in{\mathcal{A}}}{{ \sup}}\text{ }\underset{\Pbb\in\Pset}{{\rm inf}}\Ebb^{\Pbb}\left[U(X^{\rho}_T)\right].
\end{equation}

\vspace{0.8em}
In order to find the value function and an optimal strategy, we follow the method outlined in the previous Section for the exponential utility. Therefore, we have to construct a stochastic process $R^{\rho}$ with terminal value
$$R^{\rho}_T=U\left(x+\int^T_0 X^{\rho}_s \rho_s \frac{dS_s}{S_s}\right),$$ and which satisfies Properties \ref{H0}. Then the value function will be given by $V(x)=R_0$.

\vspace{0.8em}
Applying the utility function to the wealth process yields
\begin{equation}
-\frac{1}{\gamma}\left(X^{\rho}_t\right)^{-\gamma}=-\frac{1}{\gamma}x^{-\gamma}\exp\left(-\int^t_0\gamma\rho_s dB_s-\int^t_0\gamma\rho_sb_s ds+\frac{1}{2}\int^t_0\gamma \left|\widehat{a}_s^{1/2}\rho_s\right|^2 ds\right).
\end{equation}

\vspace{0.8em}
This equation suggests the following choice
\begin{equation*}
R^{\rho}_t=-\frac{1}{\gamma}x^{-\gamma}\exp\left(-\int^t_0\gamma\rho_s dB_s-\int^t_0\gamma\rho_sb_s ds+\frac{1}{2}\int^t_0\gamma \left|\widehat{a}_s^{1/2}\rho_s\right|^2 ds+Y_t\right),
\end{equation*}
where $(Y,Z)\in \mathbb D^{\infty}_H\times\mathbb H^{2}_H $ is the unique solution of the following 2BSDE
\begin{equation}
Y_t=0-\int^T_tZ_sdBs-\int^T_t\widehat F_s(Z_s)ds+K_T-K_t, \ t\in[0,T], \ \Pset-q.s.
\end{equation}

\vspace{0.8em}
In order to get $\rm{(iii)}$ of Properties \ref{H0} for $R^{\rho} $, we have to construct $\widehat F_t(z)$ such that, for $t\in[0,T] $
\begin{equation}
\gamma\rho_tb_t-\frac{1}{2}\gamma \left|\widehat{a}_t^{1/2}\rho_t\right|^2 -\widehat F_t(Z_t)\leq -\frac{1}{2} \left|\widehat{a}_t^{1/2} (\gamma\rho_t-Z_t)\right|^2\ \text{for all } \rho\in\mathcal{A},
\end{equation}
with equality for some $\rho^*\in\mathcal{A} $. It is equivalent to
\begin{align*}
\widehat F_t(Z_t)\geq -\frac{1}{2} \gamma\left(1+\gamma\right)\left|\widehat{a}_t^{1/2}\rho_t-\frac{1}{1+\gamma}\left(-\widehat{a}_t^{1/2}Z_t+\widehat{\theta}_t\right)\right|^2-\frac{1}{2}\frac{\gamma\left|-\widehat{a}_t^{1/2}Z_t+\widehat{\theta}_t\right|^2}{1+\gamma}+\frac{1}{2}\left|\widehat{a}_t^{1/2}Z_t\right|^2,
\end{align*}
with $\widehat \theta_t:=\widehat a^{-1/2}_tb_{t}$.

\vspace{0.8em}
Hence, the appropriate choice for $\widehat F$ is
\begin{align}\label{ranma}
\widehat F_t(z)&=-\frac{\gamma(1+\gamma)}{2}\text{dist}^2\left(\frac{-\widehat{a}_t^{1/2}z+\widehat{\theta}_t}{1+\gamma},A_{\widehat a_t} \right)+\frac{\gamma\left|-\widehat{a}_t^{1/2}z+\widehat{\theta}_t\right|^2}{2(1+\gamma)}+\frac{1}{2}\left|\widehat{a}_t^{1/2}z \right|^2,
\end{align}

\vspace{0.8em}
and a candidate for the optimal strategy must satisfy
$$\widehat{a}_t^{1/2}\rho^*_t\in \Pi_{A_{\widehat a_t}}\left(\frac{1}{1+\gamma}\left(-\widehat{a}_t^{1/2}Z_t+\widehat{\theta}_t\right)\right),\ t\in[0,T].$$

\vspace{0.8em}
The above results are summarised in the following Theorem.

\begin{Theorem}\label{theorempower}
Assume either that the drift $b$ verifies that $\underset{0\leq t\leq T}{\sup}\No{b_t}_{\mathbb L^\infty_H}$ is small and that the set $A$ contains $0$, or that the set $A$ is $C^2$ (in the sense that its border is a $C^2$ Jordan arc). Then, the value function of the optimization problem \reff{robustpower} is given by
\begin{equation*}
V(x)=-\frac{1}{\gamma}x^{-\gamma}{\rm exp}(Y_0)\ \ {\rm for} \  x>0,
\end{equation*}
where $Y_0$ is defined as the initial value of the unique solution $(Y,Z)\in \mathbb D^{\infty}_H\times\mathbb H^{2}_H $ of the quadratic 2BSDE
\begin{equation}
\label{2BSDE1}
Y_t=0-\int^T_tZ_sdBs-\int^T_t\widehat F_s(Z_s)ds+K_T-K_t, \ t\in[0,T], \ \Pset-q.s.,
\end{equation}
where $\widehat F$ is given by \reff{ranma}.

\vspace{0.5em}
Moreover, there is an optimal trading strategy $\rho^*\in\Ac $ with the property
\begin{equation}
\label{optim2}
\widehat{a}_t^{1/2}\rho^*_t\in \Pi_{A_{\widehat a_t}}\left(\frac{1}{1+\gamma}\left(-\widehat{a}_t^{1/2}Z_t+\widehat{\theta}_t\right)\right),\ t\in[0,T],
\end{equation}
with $\widehat \theta_t:=\widehat a^{-1/2}_tb_{t}$.
\end{Theorem}

\vspace{0.5em}
\proof
The proof is very similar to the case of robust exponential utility. First it can be shown with the same arguments, that the generator $\widehat F$ satisfies the conditions of Assumption \ref{assump.hquad} or Assumption \ref{assump.hh}. Hence there exists a unique solution to the 2BSDE \reff{2BSDE1}.

\vspace{0.8em}
Let then $\rho^* $ denote the progressively measurable process, constructed with a measurable selection theorem, which realizes the distance in the definition of $\widehat F$. The same arguments as in the case of robust exponential utility show that $\rho^*\in\Ac$.

\vspace{0.8em}
Then with the choice made for $\widehat F$, we have the following multiplicative decomposition
\begin{equation*}
R^{\rho}_t=-\frac{1}{\gamma}x^{-\gamma}\Ec\left(-\int_0^t\left(\gamma\rho_s-Z_s \right)dB_s\right)e^{-\gamma K_t^\mathbb P}{\rm exp}\left(-\int^t_0 v_sds\right),
\end{equation*}
where $$ v_t=\gamma\rho_tb_t-\frac{1}{2}\gamma \left|\widehat{a}_t^{1/2}\rho_t\right|^2 -\widehat F_t(Z_t)+\frac{1}{2} \left|\widehat{a}_t^{1/2} (\gamma\rho_t-Z_t)\right|^2\leq 0,\ dt\otimes \Pbb{\rm -a.e.}$$

\vspace{0.8em}
Then since the stochastic integral $\int_0^t(\rho_s-Z_s)dB_s$ is a $\mathbb B\rm{MO}(\mathcal P_H)$ martingale, the stochastic exponential above is a uniformly integrable martingale. By exactly the same arguments as before, we have
\begin{equation*}
\underset{\Pbb'\in\mathcal P_H(s+,\mathbb P)}{\einf^{\Pbb}}\Ebb^{\Pbb'}_s\left[R^{\rho}_{t}\right]\leq R^{\rho}_{s},\ \ s\leq t,\ \Pbb-a.s.,
\end{equation*}
with equality for $\rho^*$.

\vspace{0.8em}
Hence, the terminal value $R^{\rho}_T $ is the utility of the terminal wealth of the trading strategy $\rho$. Consequently,
\begin{equation*}
\underset{\mathbb P\in\mathcal P_H}{\inf}\Ebb^{\Pbb}\left[U\left(X^{\rho}_T\right)\right]\leq R_0=-\frac{1}{\gamma}x^{-\gamma}{\rm exp}(Y_0)\ \ \text{\rm for all }\rho\in\Ac.
\end{equation*}
\ep

\begin{Remark}
Of course, the min-max property of Theorem \ref{minmax} still holds.
\end{Remark}

\vspace{0.8em}
%%%%%%%%%%%%%%%%%%%%%%%%%%%%%%%%%%%%%%%%%%%%%%%%%%%%%%%%%%%%%%%%%%%%%%%%%%%%%%%%%%%%%%%%%%%%%%%%%%%%%%%%%%
%2.4 Robust log utility
\section{Robust logarithmic utility}\label{5}
In this section, we consider the logarithmic utility function
$$U(x)={\rm log} (x),\ x>0.$$

Here we use the same notion of trading strategies as in the power utility case,  $\rho=(\rho^i)_{i=1,\ldots,d} $ denotes the part of the wealth invested in stock $i$. The number of shares of stock $i$ is given by $\frac{\rho^i_tX_t}{S^i_t}$. Then the wealth process is defined as
\begin{equation}
X^{\rho}_t=x+\int^{t}_{0}{\sum^{d}_{i=1}}\frac{X^{\rho}_s {\rho}^i_{s}}{S^i_{s}} dS^i_{s}=x+\int^{t}_{0}X^{\rho}_{s}{\rho}_s\left(dB_s+b_sds\right),\ \mathcal{P}_{H}-q.s.
\end{equation}
and the initial capital $x$ is positive.

\vspace{0.8em}
The wealth process $X^{\rho}$ can be written as
$$X^{\rho}_t=x\mathcal{E}\left(\int_0^t \rho_s(dB_s+b_sds)\right), \ t\in\left[0,T\right],\ \mathcal{P}_{H}-q.s.$$

\vspace{0.8em}
In this case, the set of admissible strategies is defined as follows

\begin{Definition}
\label{admiss.log}
Let $A$ be a closed set in $\mathbb{R}^{d}$. The set of admissible trading strategies $\mathcal{A}$ consists of all $\mathbb{R}^d$-valued progressively measurable processes $\rho$ satisfying
$$\underset{\mathbb P\in\mathcal P_H}{\sup}\mathbb E^{\mathbb{P}}\left[\int^{T}_{0}|\widehat{a}^{1/2}_{t}\rho_{t} |^{2}dt\right]<\infty,$$
 and $\rho\in A, \;dt\otimes d\Pbb-a.s., \;\forall \Pbb\in\Pset $.
\end{Definition}

\vspace{0.8em}
For the logarithmic utility, we assume the agent has no liability at time $T$ ($\xi=0$). Then the optimization problem is given by
\begin{align}
\label{optimprob3}
 \nonumber V(x)&=\underset {\rho\in\mathcal{A}}{{\rm sup}}\ \underset {\Pbb\in\Pset}{{\rm inf}}\Ebb^{\Pbb}[\log(X^{\rho}_T)]\\
&=\log(x)+\underset {\rho\in\mathcal{A}}{{\rm sup}} \text{ }\underset {\Pbb\in\Pset}{{\rm inf}}\Ebb^{\Pbb}\left[\int^T_0\rho_sdB_s+\int^T_0(\rho_s b_s-\frac{1}{2}\lvert\widehat{a}^{1/2}_s\rho_s\rvert^2)ds\right].
\end{align}

\vspace{0.8em}
We have the following theorem.

\begin{Theorem}
Assume either that the drift $b$ verifies that $\underset{0\leq t\leq T}{\sup}\No{b_t}_{\mathbb L^\infty_H}$ is small and that the set $A$ contains $0$, or that the set $A$ is $C^2$ (in the sense that its border is a $C^2$ Jordan arc). Then, the value function of the optimization problem \reff{optimprob3} is given by
\begin{equation*}
V(x)=\log(x)-Y_0\ \ {\rm for} \  x>0,
\end{equation*}
where $Y_0$ is defined as the initial value of the unique solution $(Y,Z)\in \mathbb D^{\infty}_H\times\mathbb H^{2}_H$ of the quadratic 2BSDE
\begin{equation}
\label{2BSDE2}
Y_t=0-\int^T_t Z_sdB_s-\int^T_t \widehat{F}_sds+K^{\Pbb}_T-K^{\Pbb}_t,\ t\in[0,T],\ \Pbb-a.s.,\ \forall \Pbb\in\Pset.
\end{equation}

\vspace{0.5em}
The generator is defined by
%\begin{equation*}
%\widehat{F}_s=F_s(\widehat{a}_s),
%\end{equation*}
%where
%$$F_s(a)=-\frac{1}{2}\text{$\rm{dist}$}^2(\theta_s,A_a)+\frac{1}{2}\lvert\theta_s\rvert^2, \text{ for } a\in \mathbb S_d^{>0}.$$
$$\widehat F_s=-\frac{1}{2}\text{$\rm{dist}$}^2(\widehat \theta_s,A_{\widehat a})+\frac{1}{2}\lvert\widehat\theta_s\rvert^2,$$ %\text{ for } a\in \mathbb S_d^{>0}.$$
where $\widehat \theta_t:=\widehat a^{-1/2}_tb_{t}$.
\vspace{0.5em}
Moreover, there exists an optimal trading strategy $\rho^*\in\Ac$ with the property
\begin{equation}
\label{optim3}
\widehat{a}_t^{1/2}\rho_t^*\in \Pi_{A_{\widehat a_t}}\left(\widehat{\theta}_t\right).
\end{equation}
\end{Theorem}

\vspace{0.8em}
\proof
The proof is very similar to the case of exponential and power utility. First we show that there is an unique solution to the 2BSDE \reff{2BSDE2}. We then write, for $t\in[0,T] $
\begin{equation*}
R^{\rho}_t=M^{\rho}_t+N^{\rho}_t,
\end{equation*}
where
\begin{eqnarray*}
M^{\rho}_t&=&\log (x) -Y_0+\int^t_0\left(\rho_s-Z_s\right)dB_s+K^{\Pbb}_t,\\
N^{\rho}_t&=&\int^t_0\left(-\frac{1}{2}\left|\widehat{a}^{1/2}_s\rho_s-\widehat {\theta}_s \right|^2+\frac{1}{2}\left|\widehat{\theta}_s\right|^2-\widehat{F}_s\right)ds.
\end{eqnarray*}

Then, we similarly prove that $\rho^*$, which can be constructed by means of a classical measurable selection argument, is in $\mathcal{A}$. Note in particular that $\rho^* $ only depends on $\widehat{\theta},\ \widehat{a}^{1/2}$ and the closed set $A$ describing the constraints on the trading strategies.

\vspace{0.8em}
Next, due to Definition \ref{admiss.log}, the stochastic integral in $R^{\rho}$ is a martingale under each $\Pbb$ for all $\rho\in\mathcal{A}$. Moreover, $\widehat{F}$ is chosen to make the process $N^{\rho}$ non-increasing for all $\rho$ and a constant for $\rho^*$. Thus, the minimum condition of $K^{\P}$ implies that $R^{\rho}$ satisfies $\rm{(iii)}$ of Properties \ref{H0}.

\vspace{0.8em}
Furthermore, the initial value $Y_0$ of the simple 2BSDE \reff{2BSDE2} satisfies
$$Y_0=-\underset{\P\in\mathcal P_H}{{\sup}}\Ebb^{\P}\left[\int^T_0 \widehat{F}_sds\right].$$
Hence,
$$V(x)=R^{\rho^*}_{0}(x)=\log(x)+\underset{\P\in\mathcal P_H}{{\sup}}\Ebb^{\P}\left[\int^T_0 \widehat{F}_sds\right].$$
\ep

\begin{Remark}
Of course, the min-max property of Theorem \ref{minmax} still holds. Moreover, it is an easy exercise to show that the 2BSDE has a unique solution in this case given by
$$Y_t=\underset{\mathbb P^{'}\in\mathcal P_H(t^+,\mathbb P)}{\esup^\mathbb P}\mathbb E^{\mathbb P^{'}}\left[\int_t^T\frac12\left(\text{$\rm{dist}$}^2(\theta_s,A_{\widehat a_s})-\abs{\theta_s}^2\right)ds\right],\ \mathbb P-a.s., \ \text{for all $\mathbb P\in\mathcal P_H$}.$$
\end{Remark}

%%%%%%%%%%%%%%%%%%%%%%%%%%%%%%%%%%%%%%%%%%%%%%%%%%%%%%%%%%%%%%%%%%%%%%%%%%%%%%%%%%%%%%%%%%%%%%%%%%%%%%%%%%%%%%%%%%%%%%%%%%%%%%%%%%%
\section{Examples}\label{6}

In general, it is difficult to solve BSDEs and 2BSDEs explicitly. In this section, some examples with an
explicit solution will be given. In particular, we show how the optimal probability measure is chosen. In all our examples, we will work in dimension one, $d=1$.

\vspace{0.8em}
First,  robust exponential utility is dealt with. We consider the case where there are no constraints on trading strategies, that is $A =\R$. Then the associated $2$BSDE has a generator which is linear in $z$. In the first example, we consider a deterministic terminal liability $\xi$ and show that our result can be compared with the one obtained by solving the HJB equation in the standard Merton's approach, working with the probability measure associated with the constant process $\overline{a}$. In the second example, we show that with a random payoff $\xi=-B_T^2$, where $B$ is the canonical process, we end up with an optimal probability measure which is not of Bang-Bang type (Bang-Bang type means that, under this probability measure, the density of the quadratic variation $\widehat{a}$ takes only the two extreme values, $\underline{a}$ and $\overline{a}$). We emphasize that this example does not have real financial significance, but  nonetheless shows that one cannot expect
 the optimal probability measure to depend only on the two bounds for the volatility unlike with option pricing in the uncertain volatility model.

\subsection{Example 1: Deterministic payoff}\label{ex.1}
In this example, we suppose that $b$ is a constant in $\R$.
From Theorem \ref{theoremexp}, we know that the value function of the robust maximization problem %\reff{valuefunction}
is given by
$$V^{\xi}(x)=-{{\rm exp}}\left(-\beta\left(x-Y_0\right)\right),$$
where $Y$ is the solution of a 2BSDE with quadratic generator. When there are no constraints, the 2BSDE can be written as follows
\begin{equation*}
Y_{t}=\xi-\int^{T}_{t}Z_{s}dB_{s}-\int^{T}_{t}\widehat{F}_s(Z_{s})ds+K^{\Pbb}_{T}-K^{\Pbb}_{t},\ \Pbb-a.s.,\ \forall \Pbb\in\Pset.
\end{equation*}

and the generator is given by
$$\widehat F_t(z):=F_t(\omega,z,\widehat a)=bz+\frac{b^2}{2\beta \widehat a}.$$% \text{ for } a\in \mathbb S_d^{>0}.$$

\vspace{0.8em}
Then the corresponding BSDEs can be solved explicitly with the same generator under each $\P$. Let $$M_{t}= e^{-\int^{t}_{0}\frac{1}{2}b^2 \widehat{a}_s^{-1}ds-\int^{t}_{0}b \widehat{a}_s^{-1}dB_s}.$$ By applying It\^o's formula to $y^{\P}_tM_{t} $,
we have
$$y^{\P}_0=\E^{\P}\left[\xi M_T-\frac{b^2}{2\beta}\int^T_0\widehat{a}^{-1}_{s}M_sds \right].$$
Since $\underline a\leq \widehat a\leq \overline a$, we derive that
$$y^{\P}_0\leq \xi-\frac{1}{2\beta}\frac{b^2}{\overline{a}}T.$$
Therefore, by the representation of $Y$, we have
$$Y_0\leq \xi-\frac{1}{2\beta}\frac{b^2}{\overline{a}}T.$$
Moreover, under the specific probability measure $\P^{\overline a}\in \Pset$, we have $$y_0^{\P^{\overline a}}=\xi-\frac{1}{2\beta}\frac{b^2}{\overline{a}}T.$$

It implies that $Y_0=y_0^{\P^{\overline a}}$, which means that the robust utility maximization problem is degenerated and is equivalent to a standard utility maximization problem under the probability measure $\P^{\overline{a}}$. This result is discussed in more details in Example \ref{ex.3} below.
%%%%%%%%%%%%%%%%%%%%%%%%%%%%%%%%%%%%%%%%%%%%%%%%%%%%%%%%%%%%%%%%%%%%%%%%%%%%%%%%%%%%%%%%%%%%%%%%%%%%%%%%%%%%%%%%%%%%%%%%%%%%%%%%%%%%%%%%%%%
\subsection{Example 2 : Non-deterministic payoff}

In this subsection, we consider a non-deterministic payoff $\xi=-B_T^2$. As in the first example, there are no constraints on trading strategies. Then, the 2BSDE has a linear generator. We can verify that $-B_T^2$ can be written as the limit under the norm $\No{\cdot}_{\mathbb L^{2,\kappa}_H}$ of a sequence which is in $\rm{UC_b(\Omega)}$, and thus is in $\Lc^{2,\kappa}_{H}$, which is the terminal condition set for 2BSDEs with Lipschitz generator (these sets are defined in Section \ref{space}). Here, we suppose that $b$ is a deterministic continuous function of time $t$.

\vspace{0.8em}
By the same method as in the previous example, let $$M_t= e^{-\int^{t}_{0}\frac{1}{2}b_s^2 \widehat{a}_s^{-1}ds-\int^{t}_{0}b_s \widehat{a}_s^{-1}dB_s},$$
then we obtain $$y^{\P}_0=\E^{\P}\left[-M_T B_T^2-\int^T_0\frac{b_s^2}{2\beta}\widehat{a}^{-1}_{s}M_sds \right].$$

By applying It\^o's formula to $M_t B_t$, we have
$$dM_tB_t=M_tdB_t+B_tdM_t-b_tM_tdt.$$

Since $b$ is deterministic, by taking expectation under $\P$ and localizing if necessary, we obtain $$\E^{\P}\left[M_TB_T\right]=\E^{\P}\left[ -\int^T_0b_tM_tdt\right]=-\int^T_0b_tdt.$$

Again, by applying It\^o's formula to $-M_t B_t^2$, we have
$$-dM_tB_t^2=-2M_tB_tdB_t-B_t^2dM_t-\widehat{a}_tM_tdt+2b_tM_tB_tdt.$$

Therefore $y^{\P}_0$ can be rewritten as
$$y^{\P}_0=\E^{\P}\left[\int^T_0-M_t\left(\widehat{a}_{t} + \frac{b_t^2}{2\beta\widehat{a}_t}\right)dt \right]-\int^T_02b_t\left(\int^t_0b_sds\right)dt.$$

\vspace{0.8em}
By analysing the map $g:x\in\R^+\longmapsto x-\frac{b^2_t}{2\beta x}$, we know that $g'(x)=1-\frac{b^2_t}{2\beta x^2}$, implying that $g$ is non-decreasing when $x^2 \geq \frac{b^2_t}{2\beta}$.

\vspace{0.8em}
Let it now be assumed that $b$ is a deterministic positive continuous and non-decreasing function of time $t$ such that $$\frac{b_0^2}{2\beta}\leq \underline{a}^2\leq \overline{a}^2\leq\frac{b_T^2}{2\beta}.$$

\vspace{0.8em}
Let $\underline{t}$ be such that $\frac{b^2_{\underline{t}}}{2\beta}=\underline{a}$ and $\overline{t}$ be such that $\frac{b^2_{\overline{t}}}{2\beta}=\overline{a}$, and define $$a^*_t:=\underline{a}\textbf{1}_{0\leq t< \underline{t}}+\frac{b_t}{\sqrt{2\beta}} \textbf{1}_{\underline{t}\leq t< \overline{t}}+\overline{a}\textbf{1}_{\overline{t}\leq t\leq T},\ 0\leq t\leq T,$$ then as in Example \ref{ex.1}, we can show that $\P^{a^*}$ is an optimal probability measure, which is not of Bang-Bang type.

%%%%%%%%%%%%%%%%%%%%%%%%%%%%%%%%%%%%%%%%%%%%%%%%%%%%%%%%%%%%%%%%%%%%%%%%%%%%%%%%%%%%%%%%%%%%%%%%%%%%%%%%%%%%%%%%%%%%%%%%%%%%%%%%%%%%%%%%%%%%%%
\subsection{Example 3 : Merton's approach for robust power utility }\label{ex.3}

Here, we deal with robust power utility. As in Example \ref{ex.1}, we suppose that $b$ is a constant in $\R$ and $\xi=0$. First, we consider the case where $A=\R$. From Theorem \ref{theorempower}, $\widehat F_t(z)$ can be rewritten as $$\widehat F_t(z)=\frac{\gamma\left|-\widehat{a}_t^{1/2}z+b\widehat{a}_t^{-1/2}\right|^2}{2(1+\gamma)}+\frac{1}{2}\left|\widehat{a}_t^{1/2}z \right|^2,$$
which is quadratic and linear in $z$.

\vspace{0.8em}
Then the corresponding BSDEs can be solved explicitly under each probability measure $\P$. We use an exponential transformation and let $$\alpha:=1+\frac{\gamma}{1+\gamma},\ y'^{\P} :=e^{-\alpha y^{\P}},\  z'^{\P} :=e^{-\alpha y^{\P}} z^{\P}.$$

By applying It\^o's formula, we know that $(y'^{\P},z'^{\P})$ is the solution of the following linear BSDE
$$dy'^{\P}_t=-\alpha y'^{\P}_t\left[\frac{\gamma}{2(1+\gamma)}\left(b^2\widehat {a}_t^{-1}-2b z_t^{\P}\right)dt+z'^{\P}_tdB_t\right],$$
with the terminal condition $y'^{\P}_T=1$.

\vspace{0.8em}
For $t\in[0,T],$ let $$\lambda_t:=\frac{\alpha\gamma}{2(1+\gamma)}b^2 \widehat a_t^{-1},\ \eta_t:= -\frac{\gamma}{2(1+\gamma)}2b \widehat a_t^{-1/2},\ {\rm and} \ M_t:= e^{\int^t_0\lambda_s-\frac{\eta_s^2}{2}ds+\int^t_0\widehat{a}_s^{-1/2}\eta_sdB_s}.$$

\vspace{0.8em}
By applying It\^o's formula to $y'^{\P}_tM_t $, we obtain
$$y'^{\P}_t=\E^{\P}_{t}\left[M_T/M_t \right],\text{ so }y^{\P}_0=-\frac{1}{\alpha}\ln\left( \E^{\P}\left[M_T \right]\right).$$

Since $\underline a\leq \widehat a\leq \overline a$, we derive that $$y^{\P}_0\leq -\frac{\gamma}{2(1+\gamma)}\frac{b^2}{\overline{a}}T.$$
Thus by the representation of $Y$, we have $$Y_0\leq -\frac{\gamma}{2(1+\gamma)}\frac{b^2}{\overline{a}}T.$$
Moreover, under the specific probability measure $\P^{\overline a}\in \Pset$, we have $$y_0^{\P^{\overline a}}=-\frac{\gamma}{2(1+\gamma)}\frac{b^2}{\overline{a}}T.$$
It implies that $Y_0=y^{\P^{\overline a}}_0$. Thus, the value of the robust power utility maximization problem is $$V(x)=-\frac{1}{\gamma}x^{-\gamma}\exp\left(Y_0\right).$$

As in Example \ref{ex.1}, the robust utility maximization problem degenerates, and becomes a standard utility maximization problem under the probability measure $\P^{\overline{a}}$. In order to shed more light on this somehow surprising result, we first recall the HJB equation obtained by Merton \cite{mer69} in the standard utility maximization problem
$$-\frac{\partial v}{\partial t}-\underset{\delta\in A}{\sup}\left[\Lc^{\delta,\alpha} v(t,x) \right]=0,$$
together with the terminal condition
$$v(T,x)=U(x):=-\frac{x^{-\gamma}}{\gamma},\ x\in\R_+ ,\ \gamma>0,$$
where $$\Lc^{\delta,\alpha} v(t,x)=x\delta b\frac{\partial v}{\partial x} + \frac{1}{2}x^2\delta^2\alpha\frac{\partial^2 v}{\partial x^2},$$ with a constant volatility $\alpha^{1/2}$.

\vspace{0.8em}
%Finally, notice that
It turns out that, when $A = \R$, the value function is given by $$v(t,x)=\exp\left(\frac{b^2}{2\alpha}\frac{-\gamma}{(1+\gamma)}(T-t)\right)U(x),\ (t,x)\in\left[0,T\right]\times\R_+ .$$

\vspace{0.8em}
Let $\alpha=\overline{a}$, we have $v(0,x)=V(x)$, which is the result given by our 2BSDE method. Intuitively and formally speaking (in the case of controls taking values in compact sets, it has actually been proved under other technical conditions in \cite{tz} that the solution to the stochastic game we consider is indeed a viscosity solution of the equation below, see also Remark \ref{tev}), the HJB equation for the robust maximization problem should then be

$$-\frac{\partial v}{\partial t}-\underset{\delta\in A}{\sup}\text{ }\underset{\alpha\in\left[\underline{a},\overline{a}\right] }{\inf}\left[\Lc^{\delta,\alpha} v(t,x) \right]=0$$
together with the terminal condition
$v(T,x)=U(x),\ x\in\R_+ $.

\vspace{0.8em}
Note that the value function obtained from our $2$BSDE approach solves the above PDE, confirming the intuition that it is the correct PDE to consider in this context. Now assume that $A=\R$. If the second derivative of $v$ is positive, then the term
$$\underset{\delta\in A}{\sup}\ \underset{\alpha\in[\underline a,\overline a]}{\inf}\left[\Lc^{\delta,\alpha} v(t,x) \right],$$
becomes infinite, so the above PDE has no meaning. It implies that $v$ should be concave. Then $\overline{a}$ is the minimizer. It explains why the robust utility maximization problem degenerates in the case $A=\R$. From a financial point of view, this is the same type of result as in the problem of superreplication of an option with convex payoff under volatility uncertainty. Then, similarly to the so-called robustness of the Black-Scholes formula, this leads to the fact that the probability measure with the highest volatility corresponds to the worst-case for the investor. However, it is clear that when, for instance, we impose no short-sale and no large sales constraints (that is to say $A$ is a segment), the problem should not degenerate and the optimal probability measure switches between the two bounds $\underline{a}$ and $\overline{a}$.

\vspace{0.8em}
Finally, notice that using the language of $G$-expectation introduced by Peng in \cite{peng}, if we let $$G(\Gamma)= \frac{1}{2}\underset{\underline{a}\leq\alpha\leq\overline{a} }{\sup}\;\alpha\Gamma=\frac{1}{2}\left(\overline{a}\left(\Gamma\right)^+ - \underline{a}\left(\Gamma\right)^-\right),$$ then the above PDE can be rewritten as follows
\begin{equation}
-\frac{\partial v}{\partial t}+\underset{\delta\in A}{\inf}\left[\Lc^{\delta,\underline{a},\overline{a}} v(t,x) \right]=0,
\label{pde}
\end{equation}
where
$$\Lc^{\delta,\underline{a},\overline{a}} v(t,x)= x^2\delta^2G\left(-\frac{\partial^2 v}{\partial x^2}\right).$$

Then, our PDE plays the same role for Merton's PDE as the Black-Scholes-Barenblatt PDE plays for the usual Black-Scholes PDE, by replacing the second order derivative terms by their non-linear versions.

\vspace{0.8em}
\begin{Remark}
It could be interesting to consider more general constraints for the volatility process. For instance, we may hope to consider cases where $\underline a$ can become $0$ and $\overline a$ can become $+\infty$. From the point of view of existence and uniqueness of the 2BSDEs with quadratic growth considered here,  all the results still hold, since there is no uniform bound on $\widehat a$ for the set of probability measures considered in \cite{posmz1} (see Definition $2.2$). However, the boundedness assumption is crucial to retain the BMO integrability of the optimal strategy and thus also crucial for our proofs. We think that without it, the problem could still be solved but by now using the dynamic programming and PDE approach that we mentioned. However, delicate problems would arise in the sense that on the one hand, if $\underline a=0$, then the PDE will become degenerate and one should then have to consider solutions in the viscosity sense, and on the other hand, if $\overline a=+\infty$, the PDE will have to be understood in the sense of boundary layers.

\vspace{0.5em}
Another possible generalization would be to consider time-dependent or stochastic uncertainty sets for the volatility. It would be possible if we were able to weaken Assumption \ref{assump.hquad}$\rm{(i)}$, which was already crucial in the proofs of existence and uniqueness in \cite{stz}. One first step in this direction has been taken by Nutz in \cite{nutz} where he defines a notion of G-expectation (which roughly corresponds to a 2BSDE with a generator equal to $0$) with a stochastic domain of volatility uncertainty.
\end{Remark}

\vspace{0.8em}
\begin{Remark}\label{tev}
In \cite{tev2}, a similar problem of robust utility maximization is considered. They consider a financial market consisting of a riskless asset, a risky asset with unknown drift and volatility and a nontradable asset with known coefficients. Their aim is to solve the robust utility maximization problem without terminal liability and without constraints for exponential and power utilities, by means of the dynamic programming approach already used in \cite{tz}. They managed to show that the value function of their problem solves a PDE similar to \reff{pde}, and also that (see Proposition $2.2$) the optimal probability measure was of Bang-Bang type, thus confirming our intuition in their particular framework. Besides, they give some semi-explicit characterization of the optimal strategies and of the optimal probability measures. From a technical point of view, the main difference between our two approaches, beyond the methodology used, is that their set of generalized controls (that is to say their set of probability measures) is compact for the weak topology, because it corresponds to the larger set $\overline{\mathcal P}_W$ defined in Section \ref{section.1rob}. It is also the framework adopted in \cite{dk}. However, as shown in \cite{denism} for instance, our smaller set $\mathcal P_H$ is only relatively compact for the weak topology. Nonetheless, working with this smaller set has no effect from the point of view of applications, and more importantly makes it possible to obtain results which are not attainable by their PDE methods, for instance with non-Markovian terminal liability $\xi$ and also when the set of trading strategies is constrained in an arbitrary closed set.
\end{Remark}

\paragraph{Acknowledgments}
The authors would like to thank two anonymous referees and an associate editor for their helpful remarks and comments.

%%%%%%%%%%%%%%%%%%%%%%%%%%%%%%%%%%%%%%%%%%%%%%%%%%%%%%%%%%%%%%%%%%%%%%%%%%%%%%%%%%%%%%%%%%%%%%%%%%%%%%%%%%%%%%%%%%%%%%%%%%%%%%%%%%%%%


\begin{thebibliography}{aa12}


\bibitem{ahs} Anderson E., Hansen L.P.,  Sargent T. (2003). A quartet of semigroups for model specification, robustness, prices of risk, and model detection. \textit{Journal of the European Economic Association}, 1:68--123.

\bibitem{alp} Avellaneda, M., Levy, A., Paras, A. (1995). Pricing and hedging derivative securities in markets with uncertain
volatilities. {\sl Applied Mathematical Finance}, 2:73--88.



\bibitem{bis} Bismut, J.M. (1973).
Conjugate convex functions in optimal stochastic control, {\sl J. Math. Anal. Appl.}, 44:384--404.

 \bibitem{elkarbar} Barrieu, P., El Karoui, N. (2011). Monotone stability of quadratic
semimartingales with applications to unbounded general quadratic BSDEs, preprint, {\it arXiv:1101.5282}. %preprint. http://arxiv.org/abs/1101.5282

\bibitem{bh} Briand, Ph., Hu, Y. (2006).
BSDE with quadratic growth and unbounded terminal value, {\sl Probab. Theory Relat. Fields}, 136:604--618.

\bibitem{bms} Bordigoni, G., Matoussi, A., Schweizer, M. (2007). A stochastic control approach to a robust utility maximization problem, {\sl "Stochastic analysis and applications", Abel Symposium}, 2:125--151.

%\bibitem{cst} \c{C}etin, U., Soner, H.M., and Touzi, N. (2007).
%Option hedging under liquidity costs, {\sl Finance and Stochastics}, 14:317--341.

\bibitem{cstv} Cheridito, P., Soner, H.M.,
Touzi, N., and Victoir, N. (2007). Second order backward stochastic
differential equations and fully non-linear parabolic PDEs, {\sl Communications on Pure and Applied Math.}, 60(7):1081--1110.

\bibitem{cvi}
Cvitani\'c, J., Karatzas, I. (1992). Convex duality in constrained portfolio optimization, {\sl Ann. Appl. Proba.}, 2:767--818.

%\bibitem{del} Dellacherie, C. (1972). Capacit\'es et processus stochastiques, {\sl Springer Verlag}.

%\bibitem{dm} Dellacherie, C., Meyer, P.-A. (1975).
%Probabilit{\'e}s et potentiel, Chapitre I {\`a} IV, {\sl Hermann, Paris}.

\bibitem{denism} Denis, L., Martini, C. (2006). A theoretical framework for the pricing of contingent claims in the presence of model uncertainty, {\sl Annals of Applied Probability}, 16(2): 827--852.

\bibitem{dk} Denis, L., Kervarec, M. (2007). Utility functions and optimal investment in non-dominated models, preprint, {\it hal:00371215}.

%\bibitem{dhp} Denis, L., Hu, M., Peng, S. (2011). Function spaces and capacity related to a Sublinear Expectation: application to G-Brownian Motion Paths, {\sl Potential Analysis}, 34: 139--161.

\bibitem{elk} El Karoui, N. (1981). Les aspects probalilistes du contr{\^o}le stochastique, {\sl Ecole d'Et{\'e} de Probabilit{\'e}s de Saint-Flour IX-1979, Lecture Notes in Mathematics, Springer, Berlin}, 876:73--238.

\bibitem{elkaroui}
El Karoui, N., Peng, S. and Quenez, M.C. (1994). Backward stochastic differential equations in finance, {\sl Mathematical Finance}, 7(1):1--71.

\bibitem{ekr}
El Karoui, N., Rouge, R. (2000). Pricing via utility maximization and entropy, {\sl Mathematical Finance}, 10:259--276.

%\bibitem{ekmn}
%El Karoui, N., Matoussi, A. and Ngoupeyou, A. (2011).  Quadratic Backward Stochastic Differential Equations In jump Diffusion Model.  preprint.


\bibitem{ej} Epstein, L.G., Ji, S. (2011). Ambiguous volatility, possibility and utility in continuous time, preprint. {\it arxiv:1103.1652}.

%\bibitem{fmm}
%Faidi, W.,  Matoussi, A., and Mnif, M. (2011). Maximization of Recursive Utilities:  A Dynamic Maximum Principle Approach, {\sl SIAM J. Financial Math.}, Vol. 2, pp. 1014-1041.

%\bibitem{fol} F\"ollmer, H. (1981). Calcul d'It\^o sans probabilit\'es, {\sl Seminar on Probability XV, Lecture Notes in Math.}, 850:143--150. Springer, Berlin.

%\bibitem{ftw}
%Fahim, A., Touzi, N., and Warin, X. (2008). A probabilistic numerical scheme for fully nonlinear PDEs, {\sl Ann. Proba.}, to appear.

\bibitem{frei}
Frei, C., Mocha, M., Westray, N. (2012). BSDEs in utility maximization with BMO market price of risk, {\sl Stoch. Proc. and their App.}, 122: 2486--2519.

\bibitem{gs} Gilboa, I., Schmeidler, D. (1989). Maximin expected utility with a non-unique prior. {\sl Journal of
Mathematical Economics}, 18:141--153.

\bibitem{gun}
Gundel, A. (2005). Robust utility maximization for complete and incomplete market
models, {\sl Finance and Stochastics}, 9:151--176.

\bibitem{hstw} Hansen, L.P., Sargent,T.J., Turmuhambetova, G.A., Williams, N. (2006). Robust control and model misspecification, \textit{Journal of Economic Theory}, 128:45--90.

\bibitem{him}
Hu, Y., Imkeller, P., and M{\"u}ller, M. (2005). Utility maximization in incomplete markets, {\sl Ann. Appl. Proba.}, 15(3):1691--1712.

%\bibitem{jmn1}
%Jeanblanc, M., Matoussi, A. and Ngoupeyou, A. (2010). Robust utility maximization in a discontinuous filtration, preprint.
%
%\bibitem{jmn}
%Jeanblanc, M., Matoussi, A. and Ngoupeyou, A. (2011). Indifference pricing of unbounded credit derivatives, preprint.

%\bibitem{kob}
%Kobylanski, M. (2000). Backward stochastic differential equations and partial differential equations with quadratic growth, {\sl Ann. Prob.} 28:259--276.
%
%\bibitem{ks} Kramkov, D. and Schachermayer, W. (1999). The asymptotic elasticity of utility functions and optimal
%investement in incomplete markets. {\sl The Annals of Applied Probability}, 9(3):904--950.

%\bibitem{lsm}
%Lepeltier, J. P. and San Martin, J. (1997). Backward stochastic differential equations with continuous coefficient, {\sl Statistics \& Probability Letters}, 32(5):425--430.

\bibitem{lyo} Lyons, F. (1995). Uncertain volatility and the risk-free synthesis of derivatives, {\sl Journal of Applied Finance},
2:117--133.


\bibitem{mer69}
Merton, R. (1969). Lifetime portfolio selection under uncertainty: the continuous time
case, {\sl Rev. Econ. Stat.}, 51:239--265.

\bibitem{nutz}
Nutz, M. (2011). Random G-expectations, Annals of Applied Probability, to appear, {\it arxiv:1009.2168}. %preprint.

\bibitem{pardpeng}
Pardoux, E. and Peng, S. (1990). Adapted solution of a backward stochastic differential equation, {\sl Systems Control Lett.}, 14:55--61.

\bibitem{peng} Peng, S. (2010). Nonlinear expectations and stochastic calculus under uncertainty, preprint, {\it arxiv:1002.4546}.

\bibitem{pos} Possama\"i, D. (2010).
Second order backward stochastic differential equations with continuous coefficient, preprint, {\it arXiv:1201.1049.}

\bibitem{posmz1}
Possama\"i, D., Zhou, C. (2012). Second order backward stochastic differential equations with quadratic growth,  preprint, {\it arXiv:1201.1050}.

%\bibitem{pst} Possama\"i, D., Soner, H.M., and Touzi, N. (2011).
%Large liquidity expansion of the superhedging costs, {\sl Asymptotic Analysis: Theory, Methods and Appl.}, to appear.

\bibitem{plis}
Pliska, S.R. (1986). A stochastic calculus model of continuous trading: optimal portfolios, {\sl Math. Operations Research}, 11:371--382.

%\bibitem{que} Quenez, M. (2004). Optimal portfolio in a multiple-priors model.{\sl Progress in Probability}, 58:291--321.
%
%\bibitem{sch}
%Schied, A. (2005). Optimal investments for risk- and ambiguity-averse preferences: a
%duality approach, preprint, TU Berlin.

\bibitem{sw} Schied, A., Wu, C.-T. (2005). Duality theory for optimal investments under model uncertainty, {\sl Statistics \& Decisions}, 23: 199--217.

\bibitem{ski} Skiadas, C. (2003). Robust control and recursive utility, {\sl Finance and Stochastics}, 7:475--489.

%\bibitem{st} Soner, H.M., and Touzi, N. (2007).
%The dynamic programming equation for second order stochastic target
%problems, {\sl SIAM J. on Control and Optim.}, 48(4):2344-2365.

\bibitem{stz}
Soner, H.M., Touzi, N., Zhang J. (2010). Wellposedness of second order BSDE's, {\sl Probability Theory and Related Fields}, to appear.

%\bibitem{stz2}
%Soner, H.M., Touzi, N., Zhang J. (2010). Dual formulation of second order target problems, {\sl Ann. of App. Prob.}, to appear.

\bibitem{stz3}  Soner, H.M.,   Touzi, N.,  Zhang, J. (2011). Quasi-sure stochastic analysis through aggregation. {\sl Elect. Journal of Prob.}, 16:1844--1879.

\bibitem{tz} Talay, D. and Zheng, Z. (2002). Worst case model risk management, {\sl Finance and Stochastics}, 6:517--537.

\bibitem{tev} Tevzadze, R. (2008). Solvability of backward stochastic differential equations with quadratic growth, {\sl Stoch. Proc. and their App.}, 118:503--515.

\bibitem{tev2} Tevzadze, R., Toronjadze, T., Uzunashvili, T. (2012). Robust utility maximization for diffusion market model with misspecified coefficients, {\sl Finance and Stochastics}, to appear.


\bibitem{vnm} Von Neumann, J. and Morgenstern, O. (1944). Theory of games and economic behavior. {\sl Princeton University Press}.

\bibitem{zar}
Zariphopoulou, T. (1994). Consumption-investment models with constraints, {\sl SIAM J. Control and Optimization}, 32(1):59--85.
\end{thebibliography}
\end{document}